\input amstex

\def\b1{\text{\bf 1}}

\def\CA{{\Cal A}}
\def\CB{{\Cal B}}
\def\CC{{\Cal C}}
\def\CD{{\Cal D}}

\def\CE{{\Cal E}}
\def\CH{{\Cal H}}
\def\CI{{\Cal I}}

\def\CL{{\Cal L}}
\def\CM{{\Cal M}}
\def\CN{{\Cal N}}

\def\CP{{\Cal P}}
 \def\CQ{{\Cal Q}}
 
\def\CS{{\Cal S}}
\def\CT{{\Cal T}}

\def\gr{\text{gr}}

\def\Hom{\text{Hom}}

\def\#{\,\check{}}

\def\n{{\natural}}

\def\Ker{\text{Ker}}

\def\Spec{\text{Spec}}

\def\limleft{\mathop{\vtop{\ialign{##\crcr
  \hfil\rm lim\hfil\crcr
  \noalign{\nointerlineskip}\leftarrowfill\crcr
  \noalign{\nointerlineskip}\crcr}}}}
\def\limright{\mathop{\vtop{\ialign{##\crcr
  \hfil\rm lim\hfil\crcr
  \noalign{\nointerlineskip}\rightarrowfill\crcr
  \noalign{\nointerlineskip}\crcr}}}}

\def\hra{\hookrightarrow}
\def\iso{\buildrel\sim\over\rightarrow} 

\def\lra{\longrightarrow}

\parskip=6pt

\documentstyle{amsppt}
\document
\magnification=1100
\NoBlackBoxes

\bigskip

\centerline{\bf TOPOLOGICAL $\CE$-FACTORS }

\bigskip

\centerline {A.~Beilinson}

\centerline {University of Chicago}
\bigskip

\centerline{\it To Bob MacPherson on his 60th birthday}
\medskip
 
\bigskip
\centerline {\bf Introduction} 

\medskip 

 {\bf 0.1.} A perfect complex $P$ of $R$-modules yields a homotopy point $[P]$ of the $K$-theory 
 spectrum $K(R)$.  The Euler characteristics, i.e., the class $\chi (P)$ of $P$ in $ K_0 (R)=\pi_0 
 K(R)$, is the connected component where $[P]$ lies, so $[P]$ can be considered as an ``animation" 
 of  $\chi (P)$.  When $R$ is commutative,  the determinant sends the fundamental groupoid of $K(R)$ to the groupoid $\CL (R)$ of graded super $R$-lines; this is  a morphism of the Picard groupoids (see \cite{Del3}),  so $[P]$ controls, in particular,  the determinant line $\det P$. 
 The local Riemann-Roch story,  as seen in  \cite{Del3} and \cite{Gr},   unfolds within these grounds.
 
 Now let $F$ be a perfect constructible complex of sheaves of $R$-modules  on a compact real analytic manifold $X$. 
  Then the complex of $R$-modules $R\Gamma (X,F)$ is perfect. In this article we prove two (closely related) results about  $[R\Gamma (X,F)]$:

(a)  Let $SS (F)$ be the  micro-support of $F$ (see \cite{K}, \cite{KS});  this is a conic Lagrangian subvariety of the cotangent bundle $T^* X$. 
 Suppose we have a closed subset $Y\subset X$ and a continuous 1-form $\nu$ defined on $X\smallsetminus Y$ which takes values in the complement to $SS(F)$. Then  $\nu$ can be used to  ``localize" $[R\Gamma (X,F)]$ at $Y$: namely, there is
 a natural homotopy point $\varepsilon_{\nu Y} (F)$ of $K(R)$, which is determined by the restrictions  of $F$ and $\nu $ to any neighborhood of $Y$, together with a natural identification of $\varepsilon_{\nu Y} (F)$ with $[R\Gamma (X,F)]$. If $Y$ is
 disjoint union of finitely many $Y_\alpha$'s, then $\varepsilon_{\nu Y} (F)= \Sigma\,\varepsilon_{\nu Y_\alpha} (F)$, hence  $$[R\Gamma (X,F)] =\Sigma\,\varepsilon_{\nu Y_\alpha} (F).\tag 0.1.1$$ If $R$ is commutative, then we get the {\it $\varepsilon$-lines}  $\CE_{\nu Y_\alpha} (F):= \det \varepsilon_{\nu Y_\alpha} (F)$ and  the {\it $\varepsilon$-factorization isomorphism} $$\det R\Gamma (X,F)\iso \otimes_\alpha   \, \CE (F)_{\nu Y_\alpha} . \tag 0.1.2$$

{\it Example.}  If 
 $X$ is a circle, $\nu = d\theta$, then the $\varepsilon$-factorization looks as follows. Suppose $Y= \{ y_\alpha \}$ is a finite non-empty set; we order it clockwise, so $\alpha \in \Bbb Z /n$. Consider the open intervals $I_\alpha := (y_{\alpha -1} ,y_{\alpha })$
and semi-closed ones $I^-_\alpha := (y_{\alpha -1} ,y_{\alpha }]$. Then $\CE (F )_{\nu y_\alpha}=\det R\Gamma_c (I_\alpha^- , F|_{I_\alpha^-} )$, and (0.1.2) is the composition of the standard identifications $\det R\Gamma (X,F)\iso  \det R\Gamma_c (X\smallsetminus Y , F|_{X\smallsetminus Y}) \otimes  \det \Gamma (Y, F|_Y )=\otimes_\alpha (\det R\Gamma_c (I_\alpha , F|_{I_\alpha}) \otimes \det F_{y_\alpha} )\iso \otimes_\alpha  \det R\Gamma_c (I_\alpha^- , F|_{I_\alpha^-} )$. 

(b)  When $R$ is a field,  Kashiwara defined the characteristic cycle $CC(F)$, which is an integral cycle of local nature supported on $SS(F)$, and proved that its intersection index with the zero section  $T^*_X X \subset T^* X$ equals the Euler characteristics 
 (the Dubson-Kashiwara microlocal index formula, see \cite{D}, \cite{K}, or \cite{KS} ch.~IX): $$\chi (X,F)=\langle CC(F),T^*_X X \rangle . \tag 0.1.3 $$  

We show  that  (0.1.3) admits a natural animation. Namely, there is a natural cycle $\varepsilon (F)_{SS(F)}$ with coefficients in $ K (R)$ supported on $SS(F)$,  whose intersection with $T^*_X X$ identifies canonically  with $[R\Gamma (X,F )]$. Here $R$ can be any associative algebra; if $R$ is a field, then  $\varepsilon (F)_{SS(F)}$ lifts $CC(F)$. In particular, for $R$ commutative,    we get a Dubson-Kashiwara-style description of $\det R\Gamma (X,F)$. 
 
Idea of the construction:  $[R\Gamma (X,F)]$ satisfies an additivity property: for every partition of $X$ by locally closed subanalytic pieces $\{ Z_\beta \}$ there is a natural identification  $[R\Gamma (X,F)]=\Sigma \, [R\Gamma_c (Z_\beta , F|_{Z_\beta})]$. Now $\nu$ yields a
 supply of special locally closed subsets $Z\subset X\smallsetminus Y$, called {\it lenses}, such that $R\Gamma_c (Z , F|_{Z})=0$ (see 2.4, 2.5;  in the above example, these were $I_\alpha^-$). Using them as  components of the partition, we get (a). And (b) is essentially a combination of (a) for all possible $\nu$ and $Y$.

{\it Remarks.} (i) In short, (b) says that 
 $[R\Gamma (X,F)]$ has essentially micro-local nature.

 (ii) The constructions are compatible with filtrations on $F$. In fact,  instead of treating individual $F$, in the exposition we play all the way  with the  $K$-theory spectrum of the whole category of constructible sheaves on $X$ subject to a micro-support condition.
  
  (iii)  If $\nu = df$, where $f$ is a $C^1$-function on $X$ which is constant on $Y_\alpha$, then $\varepsilon (F)_{\nu Y_\alpha}$ comes from the Morse complex at $Y_\alpha$ and (0.1.1) comes  from the  filtration on $R\Gamma (X,F)$  provided by the Morse theory. 
 Notice that the exactness of $\nu$ is essential when one wants to see the actual cohomology, and not merely the Euler characteristics, however animated (consider the example of $X=S^1$, $F$ a constant sheaf, and $\nu =d\theta$).  

(iv) One can define (0.1.2) working directly with $\CL (R)$  (and avoiding 
 $K(R)$). Nonetheless, spectra still lurk there for the Picard groupoid $\CL (R)$ is {\it not} strictly commutative.\footnote{Picard groupoids are essentially the same as spectra with $\pi_{\neq 0,1}$ vanishing, see 1.5(iv);  a Picard groupoid is strictly commutative iff  the corresponding spectrum is abelian, i.e., comes from a (length 2)  complex of abelian groups.}

{\bf 0.2.} It would be very interesting to find out if the construction of this note (purely transcendental, as it is) admits a motivic ($\ell$-adic or de Rham) counterpart. Indeed, the idea of the $\varepsilon$-factorization of $\det R\Gamma$ comes from arithmetics: for a constructible sheaf on a curve over finite field, the $\varepsilon$-factorization of the determinant of the Frobenius action on its cohomology is inherent for the Langlands reciprocity. The precise statement 
was conjectured by Deligne  \cite{Del1} and proved by Laumon \cite{L}.  To my knowledge, the geometric  $\varepsilon$-factorization, which would be an $\ell$-adic version of (0.1.2) providing the classical  $\varepsilon$-factorization by passing to the traces of the Frobenius symmetry,\footnote{Since the classical $\varepsilon$-factors also depend on an additive character of the base finite field,  their geometric counterparts should lie in an appropriate gerbe rather than be mere graded super lines.} is not properly understood\footnote{  \cite{L} gives  a construction  when $\nu = d\log f$, $f$ is a rational function,  and the base field is finite.} (to say nothing about the format of (0.1.1) or of 0.1(b)). Notice that on higher dimensional varieties in finite characteristics   the notion of micro-support for $\ell$-adic sheaves  is not developed  (but see  \cite{AS}). In the case of curves, the geometric $\varepsilon$-factorization  in the de Rham setting was constructed in \cite{Del2} (and reinvented in \cite{BBE}) using polarized determinants of certain Fredholm operators.  

 Guided by the analogy, Deligne considered in  \cite{Del2}  the problem of $\varepsilon$-factorization of the determinant of a period matrix on a curve. At the level of numbers, the result was established in \cite{BDE} by a variant of Laumon's method.  On the geometric level,  a factorization format  was envisioned  in the last expos\'e of  \cite{Del2}.  Namely, for an algebraic vector bundle with a connection  on a complex curve, one defines its Betti cohomology as the cohomology of the sheaf $F$ of tame horizontal sections (which lives on  the real blow-up $X$ of the Riemann surface at the singular points of the connection, see \cite{M}). One has  the Betti $\varepsilon$-factorization\footnote{In \cite{Del2}  the Betti $\varepsilon$-lines were defined only for  $\nu =  d\log |f|$, $f$ is a rational function.} (0.1.2).
 The period isomorphism identifies the Betti cohomology   with the de Rham one. Then Deligne's  format asks  for a canonical identification of the Betti $\varepsilon$-lines with the de Rham ones that would  factorize  the determinant of the period isomorphism.  
  
Here is another  question  proposed by Drinfeld.
As in \cite{KS}, our habitat is a smooth  variety, which does not look very natural for the story. What   intrinsic  geometry is truly relevant for the micro-local analysis of sheaves? It should make sense outside the smooth context, so that one could  play with singular spaces directly, without embedding them into smooth ones.

I am grateful to P.~Bressler, P.~Deligne,  D.~Kaledin, M.~Kapranov,  M.~Kashiwara, B.~Tsygan and, especially, to S.~Bloch, V.~Drinfeld,  and H.~Esnault  for valuable discussions,   to J.~Jardine for an e-mail correspondence, and to the referee who pointed out a better formulation of the theorem in 4.2. The research was partially supported by NSF grant DMS-0401164.

{\bf 0.3.} {\it The contents. }  In \S1 we remind a few elementary facts about spectra, and \S2 is essentially a review of some micro-support basics (with proofs)  from  \cite{KS} 5.1, 5.2 in a format
we need. Sections 3 and 4 treat, respectively,   0.1(a) and 0.1(b).

{\it Notation: } For a subset $Z$ of a topological space $X$, we denote by $\bar{Z}$ its closure, by Int$(Z)$ its interior, and by $\partial Z := \bar{Z}\smallsetminus $Int$(Z)$ the boundary. The embedding $Z\hra X$ is denoted by $i_Z = i_{ZX}$, or (if $Z$ is open) by $j_Z = j_{ZX}$.

 \bigskip

\centerline {\bf  1. A  spectral reminder}
  
\medskip

This  section is intended to alleviate a reader who, like the author, feels foreign  in the lands of homotopy theory.   This is {\it not} a review of the basics of the subject; we merely recall a few needed facts, structures, and constructions.

{\bf 1.1. }  {\it An informal homotopy comment.} A peculiar trait of the homotopy theory  is  the lack of intrinsic language. By default, people 
resort to a description of the homotopy world as a clever ``homotopy localization" of a usual  category of rigid  objects (topological spaces, simplicial sets, complexes, etc.), referred to as  model category. The latter is rather an artificial device,  an amber spyglass of the  trade,  one needs to see objects of the homotopy world. An unsettling quality of this order of things  is that  all  constructions must be performed at the rigidified level, which takes  effort and ingenuity,  and adds arbitrariness.\footnote{The search for a model category where a given construction can be performed could be a serious problem  (as the example of the symmetric monoidal structure on spectra shows).}
 
Intuitively, a homotopy world is a kind of $\infty$-category (see Remark below). Lopping off  higher homotopies (i.e., replacing the clever localization  by a stupid one), one gets a plain  category - the homotopy category of the model category. This is a desolate place where no interesting constructions can be performed. The homotopy world is a (yet unnamed) animation of the  homotopy category that hovers inbetween the latter and  highly non-canonical model categories.

 {\it Remark.}  In \cite{Gr} Grothendieck suggested  to perceive homotopy types as $\infty$-groupoids.\footnote{So the homotopy type of a topological space is its fundamental $\infty$-groupoid whose objects are points, 1-morphisms are paths between points, 2-morphisms are homotopies between paths, etc.} Unfortunately, a simple intrinsic definition of the concept is not available. 
 Arguably,  the common language of category theory may be inadequate for describing the homotopy world.\footnote{A shade of this inadequacy presents already in the plain category theory: while its force lies in the fact that one need not distinguish equivalent categories, it ostensibly asserts that for a (small) category  its set of objects is a meaningful notion.}

{\bf 1.2.} We will be interested in the homotopy world of spectra $\CS$.
For our modest purposes (we will not use the tensor structure on spectra),  the  traditional model category of spectra  from \cite{BF} will do. So  
a spectrum is a sequence $P$ of pointed simplicial sets $P_0 ,P_1, \ldots$, connected by  maps $\alpha_n : S P_n \to P_{n+1}$; here $S P:= S^1 \wedge P$ is the suspension. The category $\tilde{\CS}$ of those is closed under limits and colimits.

If $X$ is any simplicial set and $P\in\tilde{\CS}$, then the sequence $X\wedge P_i$ forms naturally a  spectrum $X\wedge P$;\footnote{Here $X\wedge P_i := X_+ \wedge P_i$ where $X_+$ is  $X$ with the marked point added. 
 Of course, $X\wedge P$ is the same as the wedge product of spectra $S^\infty X$ and $P$;  as was mentioned, the wedge products of arbitrary spectra are not needed for our modest purposes. } thus the monoidal category of simplicial sets (with $\times$ as the tensor product) acts on $\tilde{\CS}$. The endofunctor $P\mapsto X\wedge P$ admits right adjoint $Q\mapsto Q^X$. A morphism $
X\wedge P\to Q$ is referred to as an {\it $X$-family of morphisms} $P\to Q$. There is a pointed simplicial set $\CH om_{\tilde{\CS}} (P,Q)$ such that
 $$\CH om_{\tilde{\CS}} (P,Q)(X)=  \Hom_{\tilde{\CS}} (X\wedge P,Q)= \Hom_{\tilde{\CS}} (P,Q^X ),  \tag 1.2.1$$ so  $\tilde{\CS}$ is  a pointed simplicial category.
 
 {\it Example. } If $S$ is the sphere spectrum,  then $\CH om_{\tilde{\CS}}(P,Q)= Q_0$.

 For $P\in \tilde{\CS}$  its {\it homotopy groups} $\pi_i P$, $i\in \Bbb Z$, are defined as $\pi_i P := \limright \pi_{i+n} |P_n |$ where $|P_n |$ is the realization of $P_n$ (which is a pointed topological space) and the limit is taken with respect to the maps $\alpha_n$. A morphism $f:P\to Q$ is said to be a {\it weak equivalence} if it induces an isomorphism between the homotopy groups; it is a {\it cofibration} if  the maps $P_0 \to Q_0$ and $P_n \mathop\cup\limits_{SP_{n-1}} SQ_{n-1} \to Q_n$, $n\ge 1$, are all injective. The datum of weak equivalences and cofibrations defines on $\tilde{\CS}$ a structure of  {\it  simplicial model category} which is {\it stable} and {\it proper} (see \cite{GS} for a brief introduction to model categories and \cite{Hi} for a detailed exposition).

The cofibrant objects for this model structure are  $P$ with  $ SP_n \to P_{n+1}$ all injective;
the fibrant objects  are  {\it  $\Omega$-spectra}, i.e., those $Q\in\tilde{\CS}$ that each $Q_n$ is a Kan simplicial set and every map $Q_n \to \Omega Q_{n+1}$  is a weak equivalence of simplicial sets.

{\bf 1.3.} The simplicial model category structure permits to play with  spectra in much the same way as with simplicial sets. 

For example, suppose $P$ is  cofibrant,
 $Q$ is fibrant.  Then  weak equivalences between such $P$'s or such $Q$'s yield weak equivalences between the $\CH om_{\tilde{\CS}}$'s;   cofibrations between $P$'s or fibrations between $Q$'s yield Kan fibrations beween the $\CH om_{\tilde{\CS}}$'s.
  In particular, $\CH om_{\tilde{\CS}}(P,Q)$ is a Kan simplicial set. 

A corollary: for any simplicial set $X$ and a cofibrant $P$ the spectrum $X\wedge P$ is also cofibrant, and the functor $(X,P)\mapsto X\wedge P$ transforms weak equivalences between $X$'s and $P$'s to weak equivalences. Sometimes we use notation $C_\natural (X,P):= X\wedge P$; this is  the {\it homology spectrum} of $X$ with coefficients in $P$. Dually, if $Q$ is fibrant, then  $C^\natural (X,Q):=Q^X$ is fibrant, and this {\it cohomology} functor transforms weak equivalences between $X$'s and $Q$'s to weak equivalences. The homology and cohomology spectra are naturally functorial; in particular, the map $X\to$ (point) yields a canonical morphism $$tr : C_\n (X,P)\to P. \tag 1.3.1$$

 One often needs to consider the homology and cohomology with non-constant coefficients. Namely,  for a simplicial set $X$ a {\it homology type coefficient system} on $X$ is a rule $\CP$  that assigns to each $a\in X_n$ a spectrum $\CP_a \in\tilde{S}$ and to each non-decreasing map $\phi :[0,m]\to [0,n]$ a morphism $\CP_a \to \CP_{\phi (a)}$ in $\tilde{\CS}$ compatible with the composition of $\phi$'s. Then the {\it homology spectrum} $C_\n (X,\CP)$ is the colimit of a diagram formed by spectra\footnote{We denote by $\Delta_n$ the standard $n$-simplex.} $\Delta_n \wedge \CP_a$ labeled by pairs $(a,n)$ as above, and $\Delta_m \wedge \CP_a$ labeled by pairs $(a,\phi )$ as above, which are connected by the  arrows $\Delta_m \wedge \CP_{\phi (a)} \leftarrow \Delta_m \wedge \CP_a  \rightarrow \Delta_n \wedge \CP_a$ coming from $\phi$. A {\it cohomology type coefficient system} $\CQ$ on $X$ is defined in a dual manner; for such $\CQ$ one constructs its {\it cohomology spectrum} $C^\n (X,\CQ)$  as the limit of a diagram formed by spectra $\CQ_a^{\Delta_n}$ labeled by $(a,n)$ as above, and $\CQ^{\Delta_n}_a $ labeled by $(a,\phi )$ as above, which are connected by the  arrows $\CQ_a^{\Delta_n}\to \CQ_a^{\Delta_m}\leftarrow \CQ_{\phi (a)}^{\Delta_m}$ coming from $\phi$. For  constant coefficient systems $P$, $Q$, we get the previous construction. 

{\it Examples.} (i) If $X$ is a point, then $\CP$ is the same as simpicial object in the category $\tilde{\CS}$, and $\CQ$ is the same as a cosimplicial object.

(ii) Let  $\CI$ be a small category and $F: \CI \to \tilde{\CS}$ be an $\CI$-diagram of spectra. 
It yields a homology  type coefficient systems $F_\n$  on the nerve $\CN er (\CI )$.  
Namely, for  $(i_0 \to\ldots\to i_n )\in \CN er (\CI )_n$ one has $F_{\n (i_0 \to\ldots\to i_n )}:=  F(i_0 )$.
Set  hocolim$_\CI F := C_\n (\CN er (\CI ), F_\n )$.
If $F\to F'$ is a morphism of diagrams such that for each $i\in \CI$ the map $F(i)\to F'(i)$ is a weak equivalence of cofibrant spectra, then the map hocolim$_\CI F \to$ hocolim$_\CI F'$ is a weak equivalence of cofibrant spectra. Dually, $F$ yields a cohomology type coefficient system $F^\n _{(i_0 \to\ldots\to i_n )}:=  F(i_n )$,  holim$_\CI (F):= C^\n (\CN er (\CI ), F^{\n})$, and holim$_\CI$ sends pointwise weak equivalences between pointwise fibrant $F$'s to weak equivalences of fibrant spectra. These  are parts of  the story of natural model structures  on the categories of diagrams, coefficient systems, etc.

The usual constructions of the standard homotopy category, such as the cone, cylinder,  path space, and homotopy fiber space of a morphism are particular cases of the  homotopy (co)limits.

{\bf 1.4.} One gets spectra localizing $\tilde{\CS}$ with respect to weak equivalences, and the model category format  explains how to do this with minimum of pain. To produce Ho$\,\CS$, one simply localizes the category $\tilde{\CS}$ with respect to weak equivalences, which is the same as to consider the subcategory of both fibrant and cofibrant objects with $\pi_0 \CH om_{\tilde{\CS}}$ as morphisms. 
This is naturally a triangulated  category with a  non-degenerate t-structure whose core equals the category  $\CA b$ of abelian groups; the corresponding homology functor assigns to a spectrum $P$ its homotopy groups $\pi_i P$.

{\it Remark.} Among all  spectra live abelian ones, which are the same as complexes of abelian groups (see 1.5(i)); on the level of the homotopy categories, one has a faithful t-exact embedding $D (\CA b) \to $ Ho$\,\CS$ which identifies  the cores.

As a move towards the elusive clever localization $\CS$ (see 1.1),  one  considers for $P,Q\in\tilde{\CS}$, instead of true morphisms $P\to Q$,  the {\it homotopy morphisms}, which are $X$-families of morphisms with base $X$ contractible. 
An identification of a homotopy morphism $f$ with base $X$ and one $g$ with base $Y$ is a homotopy morphism $h$ with base $Z$ and maps $X\to Z\leftarrow Y$ such that $f$ and $g$ are the respective pull-backs of $h$, etc. If $P$ is the spherical spectrum, we refer to homotopy morphisms $P\to Q$ as {\it homotopy points} of $Q$. In practice, the base $X$ is provided in the course of constructions and rarely specified explicitly.
A related notion is that of {\it homotopy object} of $\tilde{\CS}$, by which we understand a simplicial functor $\CB \to \tilde{\CS}$ where $\CB$ is a non-empty small simplicial category all of whose $\CH om$ spaces are contractible.

{\it Examples.} 
(i) For any $P\in\tilde{\CS}$ its {\it fibrant resolution}  is a morphism  $i: P\to P^f$ such that $i$ is an acyclic cofibration\footnote{I.e., $i$ is a cofibration  and a weak equivalence.} and $P^f$ is fibrant. Fibrant resolutions exist, and form a single homotopy object. More precisely, for any morphism $\phi : P\to Q$ and fibrant resolutions $P^f$, $Q^f$ all morphisms $\phi^f : P^f \to Q^f$ compatible with $\phi$ form a contractible simplicial subset of $\CH om_{\tilde{\CS}}(P^f,Q^f)$, so $\phi^f$ is canonically defined as a homotopy morphism. There is a dual statement for {\it cofibrant resolutions} $P^c \to P$. 
Iterating, one get fibrant-cofibrant $P\leftarrow P^c \to (P^c )^f =:P^{fc}$ and cofibrant-fibrant $  P\to P^{f}\leftarrow (P^f )^c =:P^{fc}$ resolutions; these form single homotopy objects which are both fibrant and cofibrant, together with a canonical weak equivalence $P^{fc}\to P^{cf}$. If wanted, the resolutions can be chosen functorially.

(ii) Suppose $P$, $Q$ are both fibrant and cofibrant, and $f: P\to Q$ is a weak equivalence. Let us spell out how  one inverts $f$ as a homotopy morphism. 
A {\it right homotopy inverse} to $f$ is a pair $(g_r ,h_r )$ where $g_r$ is a morphism $g_r : Q\to P$ and $h_r$ is  a homotopy $\Delta_1 \wedge Q\to Q$ between $fg_r$ and $id_Q$. The right homotopy inverses form a contractible simplicial set. The  left homotopy inverses $(g_\ell ,h_\ell )$ are defined dually; they form a contractible simplicial set as well. Thus we have canonical homotopy morphisms $g_\ell ,g_r : Q\to P$. To identify them, one uses the contractible set of collections $(g_\ell ,h_\ell ; g_r ,h_r ; \tilde{h})$ where  $(g_\ell ,h_\ell )$, $(g_r ,h_r )$ are left and right homotopy inverses to $f$ and $\tilde{h} :\Delta_2 \wedge Q\to P$ is a map which equals $g_\ell$, $g_\ell fg_r$, and $g_r$ at the vertices, $g_\ell h_r$, $h_\ell g_r$ at the two edges. 

{\it Remark.}  If $f$ is a fibration, then its true right inverses (i.e., sections) $g$ form a contractible simplicial set, as well as the pairs $(g,h)$ where $h$ is a fiberwise homotopy between $gf$ and $id_P$. There is a dual statement if $f$ is a cofibration.

(iii) Let us show that spectra carry a canonical homotopy sum operation. Take any $P$ which is both fibrant and cofibrant,  and let $I$ be any finite set. Consider a morphism $e : \bigvee_I P\to P^I$ whose components $e_{i i'} : P\to P$ equal $id_P$ if $i=i'$ and are the trivial maps otherwise. Then $e$ is an acyclic cofibration, so
all morphisms $\Sigma_I : P^I \to P$ such that the composition $\Sigma_I e$ is the map which has $id_P$ as the components, form a contractible simplicial set $E(P)_I$. Thus we have  a canonical
single homotopy morphism $\Sigma =\Sigma_I : P^I \to P$. When $I$ changes, the  spaces $E(P)_I$   form in an evident way an operad which acts on $P$.

{\bf 1.5.} {\it Examples.} (i) Any complex of abelian groups $A$ defines an Eilenberg-Mac\,Lane spectrum $EM(A)$. Namely, for each $n\ge 0$ the complex $\tau_{\le 0} (A[n])$ can be seen, by Dold-Puppe, as a commutative simplicial group which we denote by $A_n$. There is an evident identification $A_n \iso \Omega A_{n+1}$, so $A_n$ form an $\Omega$-spectrum. Notice that for any simplicial set $X$ one has $C_\natural (X, EM(A))= EM (C_\natural (X,A))$, $C^\natural (X, EM (A))= EM(C^\natural (X,A))$, where $
C_\natural (X,A)$ and $C^\natural (X,A)$ are the usual simplicial homology and cohomology complexes.

(ii) Let $P$ be as in Example (iii) in 1.4. The operad $E_P$ acts on the (Kan) simplicial set $P_0$; the corresponding  operation on $\pi_0 (P_0 )$ is the group operation on $\pi_0 (P)$. Conversely, suppose we have an  operad $E$ with contractible terms. Then any $E$-space, i.e., 
a Kan simplicial set equipped with an action of  $E$,  such that $\pi_0$ is a group can be lifted in a natural (and  homotopically unique) way to an $\Omega$-spectrum with $\pi_{<0}=0$. 

(iii)  If $\CM$ is a small symmetric monoidal category, then  its nerve $\CN er  ( \CM )$ is naturally an $E$-space;  here  $E$ is the nerve of an operad in groupoids whose $n$th term is the contractible groupoid of all ``n-fold operation"  functors $\CM^n \to \CM$. Thus $\CN er ( \CM )$ is  a spectrum if  $\pi_0 \,\CN er (\CM )$ is a group.

(iv)  The latter happens if $\CM$  is a Picard groupoid.\footnote{Which means that each morphism in $\CM$ is invertible, and each object is invertible with respect to the operation.} The only non-trivial homotopy groups of $\CN er (\CM )$ are $\pi_0 $ (which is the group of isomorphism classes of objects of $\CM$) and $\pi_1 $ (which is the automorphism group of any object in $\CM$).  Conversely, for any $P$ as in (ii) the $E_P$-operation defines on the fundamental groupoid $\Pi P:=\Pi P_0$ the structure of a Picard groupoid; if $\pi_{<0}P=0$, then one has a natural homotopy morphism of spectra $P\to \CN er (\Pi P )$ which is an isomorphism on $\pi_0$, $\pi_1$. These constructions identify homotopically the 2-category of Picard groupoids with that of spectra having all homotopy groups but $\pi_1$, $\pi_0$  trivial.

{\bf 1.6.} For the next subject, see \cite{J1}, \cite{J2}.   Let $\CT$ be a small category.
A {\it presheaf $P$ of  spectra} on $\CT$ is the same as a $\CT^\circ$-diagram of  spectra. Such objects form a simplicial category $PS h (\CT ,\tilde{\CS})$. Our $P$ yields presheaves $\pi_i P$ of homotopy groups. A morphism of presheaves $P\to Q$ is said to be {\it objectwise weak equivalence} if the morphisms $\underline{\pi}_i P\to \underline{\pi}_i Q$ are all isomorphisms; it is a cofibration if the maps $P(U)\to Q(U)$ are cofibrations for all $U\in\CT$.  These data define on $PS h (\CT ,\tilde{\CS})$ 
a stable proper simplicial model category structure.

Suppose that  $\CT$ is equipped with a Grothendieck topology. For a presheaf $P$ of  spectra we denote by $\underline{\pi}_i P$ the sheafification of the presheaf  $\pi_i P$. A morphism of presheaves $P\to Q$ is said to be {\it local weak equivalence} if all $\underline{\pi}_i P\to \underline{\pi}_i Q$ are  isomorphisms. The local weak equivalences and the  above cofibrations  define on $PS h (\CT ,\tilde{\CS})$ a stable proper simplicial model category structure (the ``sheafified" model structure). Its fibrant objects can be characterized by a homotopy descent property with respect to all hypercoverings in $\CT$ (see  \cite{DHI}). The corresponding homotopy category is naturally a t-category whose core is the category of sheaves of abelian groups on $\CT$ and the cohomology functor  $P\mapsto \underline{\pi}_{-i} P$.

The embedding of the presheaves of abelian spectra, i.e., 
 complexes of presheaves of abelian groups, into $\CP\CS h (\CT ,\tilde{\CS})$ yields a t-exact functor from $D(\CT ) $ to the homotopy  category of $\CP\CS h (\CT ,\tilde{\CS})$ with respect to the ``sheafified" model structure which induces an equivalence between the cores. Here $D(\CT )$ is the derived category of sheaves of abelian groups on $X$ (which is the same as the localization of the homotopy category of complexes of presheaves modulo local quasi-isomorphisms).

For a presheaf $P$ consider its fibrant resolution $P\to P^f$ for the ``sheafified"  model structure (it is uniquely defined as a homotopy object and can be chosen functorially, see Example (i) in 1.4), and define $R\Gamma (\CT ,P):= $ lim$_{\CT^\circ}P^f \iso $  holim$_{\CT^\circ} P^f$. 

{\it Remarks.} (i) $R\Gamma$ yields a triangulated functor between the homotopy categories.
For presheaves of abelian spectra, $R\Gamma$ is the usual cohomology functor.

(ii) Let $\CT' \subset \CT$ be a subcategory such that every object of $\CT$ admits a covering by objects of $\CT'$, and every $\CT$-covering of an object of $\CT'$ admits a refinement which lies in $\CT'$. Equip $\CT'$ with induced Grothendieck topology; then the restriction functor  
$\CP\CS h (\CT ,\tilde{\CS})\to \CP\CS h (\CT' ,\tilde{\CS})$ is compatible with the model categories structures, induces equivalences between the localized categories, and for any $P\in \CP\CS h (\CT ,\tilde{\CS})$ one has $R\Gamma (\CT ,P)\iso R\Gamma (\CT' ,P|_{\CT'})$.

{\bf 1.7.}  Let $X$ be a locally compact topological space. Any cofibrant  spectrum
 $P$ yields  a cofibrant  presheaf $P^!_X$  on $X$,   $U\mapsto P_X^! (U)= C_\n (X,X\smallsetminus U;P):=\CC one (C_\n (X\smallsetminus U ,P)\to C_\n (X,P))$, the restriction maps are the evident ones.\footnote{To fit into simplicial format, we replace each topological space by a natural simplicial approximation (say, all continuous simplices).} This is a spectral version of the dualizing complex with coefficients in $P$. One has the corresponding Poincar\'e duality theorem:

\proclaim{\quad Proposition} Suppose $X$ is compact and  has finite homological dimension. Then the morphism $C_\n (X,P)\to R\Gamma (X,P^!_X )$ is a weak equivalence.
 \endproclaim

 {\it Proof.}  The fact is standard if $P$ is abelian. The case of general $P$ can be reduced to the abelian one as follows. The canonical filtration on $P$ defined by the t-structure (i.e., the Postnikov tower) yields  filtrations on  $C_\n (X,P)$, $ R\Gamma (X,P^!_X )$, and our morphism is compatible with them. The corresponding spectral sequences with $E^{p,q}_1 = \pi_{-q} C_\n (X, EM( \pi_{-p} P))$ and $ E^{p,q}_1 = \pi_{-q} R\Gamma (X, EM(\pi_{-p} P))$ converge, respectively, to $\pi_{-p-q} C_\n (X,P)$ and $\pi_{-p-q}R\Gamma (X,P^!_X )$ due to the finiteness condition. Our morphism yields an isomorphism on $ E^{p,q}_1$, and we are done. 
 \hfill$\square$
 
{\bf 1.8. } For  details on $K$-theory spectra, see \cite{DS}, \cite{C} and references therein.

 The basic construction of $K$-theory assigns to every small pretriangulated DG category $\CA$ a spectrum $K(\CA )$, to each DG functor $\phi :\CA \to\CB$ a homotopy morphism of spectra $K(\phi ): K(\CA )\to K(\CB )$,  to any $d$-morphism\footnote{A $d$-morphism between objects of a DG category, or between DG functors, is a morphism of degree 0 which commutes with the differential.}  of DG functors $\phi \to \phi'$ such that for every $M\in\CA$ the morphism $\phi (M)\to \phi' (M)$ is a homotopy equivalence,  a homotopy  between $K(\phi )$ and $K(\phi' )$, etc. We can (and will) assume that $K(\CA)$ is fibrant and cofibrant. Set $K_i (\CA ):= \pi_i K(\CA )$. 

 In fact, $K(\CA )$ is constructed as an $E_\CA$-space (see 1.5(ii)) where $E_\CA$ is the nerve of an operad in groupoids whose $n$th term is the contractible groupoid of  all direct sum DG functors $\CA^n \to \CA$. We  do not need a precise construction of  $K (\CA )$ (consider it as a black box). The next properties will be of use:

(i)  The $K$-functor commutes with finite direct products (and the  above operad $E_\CA$ acts on $K(\CA )$ via the direct sum functors $K(\CA )^n =K(\CA^n )\to K(\CA)$).

(ii) $K_0 (\CA )$ is the Grothendieck group of the  triangulated category $\CA^{tri}$ defined by $\CA$. 

(iii) If $\phi^{tri} :\CA^{tri}\to \CB^{tri}$ is an equivalence of categories, then $K(\phi )$ is a weak equivalence. In particular, as a homotopy object, $K(\CA )$ makes sense whenever $\CA^{tri}$ is essentially small. By abuse of notation, we can write $K(\CA )=K(\CA^{tri})$.

(iv) For $n\ge 1$ let   $\CA^{(n)}$ be the DG category of $n$-filtered objects, i.e., diagrams $A_\cdot = (A_1 \to \ldots \to A_n )$ in $\CA$ (here $\to$ are $d$-morphisms, i.e., closed morphisms of degree 0 in $\CA$). For $i=1,\ldots ,n$ set $\gr_i A_\cdot := \CC one (A_{i-1}\to A_i )$.\footnote{For $i=1$ we set $A_{-1}:=0$.}
Then the DG functor $\gr_\cdot =(\gr_i ) : \CA^{(n)}\to \CA^n$ yields a weak equivalence  $K( \CA^{(n)})\to K(\CA )^n$. Notice that a right inverse to $\gr_\cdot $ is provided by a functor $\CA^n \to \CA^{(n)}$, $(A_i )\mapsto (A_1 \hra A_1 \oplus A_2 \hra\ldots A_1 \oplus\ldots\oplus A_n )$. 

This property will be used in the following form: Suppose one has a sequence $\phi$ of DG functors and $d$-morphisms $\phi_1 \to\ldots\to\phi_n$. Then the morphisms of the $K$-spectra $K(\phi_n )$ and $\sum K(\gr_i \phi )$ are naturally identified as homotopy morphisms; here $\gr_i \phi := \CC one (\phi_{i-1}\to \phi_i )$.

(v) Consider the category $\CA^\times$ whose objects are the same as objects of $\CA$ and morphisms are $d$-morphisms in $\CA$ which are homotopically invertible (i.e., are invertible in $\CA^{tri}$); this is a symmetric monoidal category with respect to $\oplus$, so $\CN er (\CA^\times )$ is an $E$-space. Then there is a natural morphism of  $E$-spaces $\CN er (\CA^\times )\to K(\CA )$. In particular, each object $A\in \CA$ yields a point in $K(\CA )$ which we denote by $[A]$.

(vi) Let $\CB^{tri} \subset \CA^{tri}$ be a thick subcategory,  $\CB \subset \CA$ the corresponding DG subcategory, and $\CA /\CB$ be the DG quotient (see \cite{Dr}), so $(\CA /\CB )^{tri}$ is the Verdier quotient $\CA^{tri}/\CB^{tri}$. Then $K(\CA^{tri} /\CB^{tri} )\iso \CC one (K(\CB^{tri} )\to K(\CA^{tri} ))$. 

If $R$ is an associative algebra,\footnote{Or, more generally, a DG algebra, or even a small DG category as in \cite{Dr}.} then $K(R)$ denotes the $K$-spectrum of the triangulated category of perfect complexes of $R$-modules (localized by quasi-isomorphisms).

(vii) Let $R$ be a commutative algebra. Denote by $\CL (R)$ the Picard groupoid of $\Bbb Z$-graded super $R$-lines. Thus objects of $\CL (R)$ are invertible $\Bbb Z$-graded $R$-modules. Locally on $\Spec\, R$ any $L\in \CL (R)$ lies in a fixed degree $\deg_L$. So $\deg_L$ is a locally constant $\Bbb Z$-valued function on $\Spec\, R$, and our $L$ is the same as a pair $(L_{o} ,\deg_L )$ where $L_{o}$ is an invertible (non-graded) $R$-module and $\deg_L$ is any function as above. By abuse of notation, we usually write $L_{o}=L$. The operation in  $\CL (R)$ is the tensor product of $\Bbb Z$-graded modules, so $\deg_{L\otimes L'}=\deg_L +\deg_{L'}$, $(L\otimes L')_{o}= L_{o}\otimes_R L'_{o}$, and the commutativity constraint for  $L\otimes L'$ equals the one for $L_o \otimes_R L'_o$ multiplied by the super sign $(-1)^{\deg_L \deg_{L'}}$.

Then the determinant line construction provides  a canonical ``character sheaf" $\lambda$ on $K(R)$, which is a local system of graded super lines together with an identification $\cdot_{xy} : \lambda_x \otimes\lambda_y \iso \lambda_{x+y}$ compatible with the associativity and commutativity constraints. Equivalently, $\lambda$ can be seen as a morphism  of the Picard   groupoids $\Pi K(R)\to \CL (R)$ (see 1.5(iv)). For a perfect $R$-complex $P$ the graded super line $\lambda_{[P]}$ equals  $\det (P)$, for a quasi-isomorphism $\alpha : P\to P'$ the monodromy of $\lambda$ along the corresponding arc in $K(R)$ equals $\det (\alpha ): \det (P)\iso \det (P' )$, and $\cdot_{[P][P']}$ is the standard isomorphism $\det (P)\otimes\det (P' )\iso \det (P\oplus P')$.

\bigskip

\centerline{\bf 2.  Morse-theoretic preliminaries}

\medskip

This section is mostly an exposition  of the story of  \cite{KS} 5.1--5.2 in a format 
convenient for our purposes. The principal fact is the proposition in 2.5 (which is a variant of 
\cite{KS} 5.1.1, 5.2.1) that locates the  micro-support of a sheaf in terms of  cohomology vanishing for certain displays of locally closed subsets of $X$, called lenses, 
controlled by conical domains in the tangent bundle. 

{\bf 2.1.}  For this subsection, $X$ is any locally compact topological space of finite cohomological dimension, and $D(X)=D(X,\Bbb Z )$ is the derived category of complexes of sheaves of abelian groups on $X$. 

Let $U,V $ be open subsets of $X$ such that $U\cap V=\emptyset$. For any sheaf $F$ the composition of the canonical maps $j_{U!} j_U^* F \to F \to j_{V*} j_V^* F$ equals 0, so $F$ carries a natural 3-step filtration $j_{U!} j_U^* F \subset \Ker (F \to j_{V*} j_V^* F)= i_{X\smallsetminus V *} i^!_{X\smallsetminus V} F \subset F$. Passing to the right derived functor, we see that every $F\in D(X)$ carries 
 a canonical 3-step {\it Morse} filtration\footnote{I.e., one has an exact functor 
 $D(X)\to DF(X)$(:= the filtered derived category) whose composition with the forgetting of the
  filtration functor $DF(X)\to D(X)$ equals $Id_{D(X)}$.} $j_{U!} j_U^* F \subset 
   i_{X\smallsetminus V *} Ri^!_{X\smallsetminus V} F \subset F$. Its top quotient 
   equals $Rj_{V*} j^*_V F$. Denote by $ M_{UV}(F)$ the middle successive quotient; 
   this is the {\it Morse complex}.
   The construction has local nature and is self-dual in the evident sense. 
 
 {\it Example.} Let $f$ be  
 a continuous $\Bbb R$-valued function on $X$.  Then each $a\in \Bbb R$ yields a pair $U= X_{f>a}$, $V=
 X_{f<a}$ as above; the corresponding Morse complex is denoted by $M_{f=a}(F)$. Here
  $X_{f<a}:=\{ x\in X : f(x)>a \}$, etc.

 We will mostly consider the case when $U= X\smallsetminus \bar{V}$;
  the Morse complex $M_V (F) := M_{UV}(F)$  is supported then on $\partial V :=\bar{V}\smallsetminus V$.

The next lemma (cf.~\cite{KS} 2.7.2) describes a situation when the vanishing of Morse complexes yields the vanishing of  global cohomology.  
 Let
$V_t$, $t\in [0,1]$, be a family of open subsets of $X$ such that $\mathop\cup\limits_{t'<t} V_{t'} = V_t$ for each $t>0$,  $\mathop\cap\limits_{t'>t} V_{t'}\subset \bar{V}_t$ for each $t<1$, and 
the subset $Z:= V_1 \smallsetminus V_0$ is relatively compact. 

\proclaim{\quad Lemma} If $M_{V_t}(F)=0$ for every $t<1$, then $R\Gamma (Z, Ri_Z^! F)=0$.
\endproclaim

{\it Proof.} We can assume that each $F^i$ is a flabby sheaf. The map $F(V_1 )\to F (V_0 )$ is surjective, and its kernel $K$ represents  $R\Gamma  (Z, Ri_Z^! F)$. 
   Take any cycle $\alpha\in K^0 $; we want to find $\beta \in K^{-1}$ such that $d\beta =\alpha$. Consider the set $\CC$ of pairs $(\beta_t ,t)$ where  $\beta_t \in F^{-1}(V_t )$ is such that $\beta_t |_{V_0} =0$ and $d\beta_t =\alpha|_{V_t}$. It is naturally ordered (one has $(\beta_t ,t)\le (\beta_{t'} ,t' )$ if $t\le t'$ and $\beta_{t'}|_{V_t}=\beta_t$) and satisfies the condition of Zorn's lemma. Thus there is a maximal element $(\beta_t ,t)\in\CC$; we want to check that $t=1$. 

If not, pick any $t' \in (t,1)$; since  $M_{V_{t'}}(F)=0$, one can find a neighborhood $U$ of $\bar{V}_{t'}$ and
$\alpha' \in F^0 (U)$ such that $d\alpha' =0$ and $\alpha'|_{V_{t'}}=\alpha|_{V_{t'}}$. 
 Since $M_{V_t}(F)=0$, one can find a neighborhood $U'$ of $\bar{V}_t$, $U'\subset U$, and $\beta' \in F^{-1}(U')$ such that $d\beta' = \alpha'|{U'}$, $\beta'|_{V_t}=\beta_t$. Since $\bar{Z}$ is compact, $U'\supset V_{t''}$ for some $t'' \in (t,t')$; set $\beta_{t''}:=\beta'|_{V_{t''}}$. Then $(\beta_{t''}, t'')\in\CC$ and $(\beta_{t''}, t'')> (\beta_t ,t)$, so
$(\beta_t ,t)$ is not maximal, q.e.d.  \hfill$\square$

{\bf 2.2.}  {\it Notation.} A closed  subset $N$ of a  finite-dimensional $\Bbb R$-vector space $T$ is said to be  {\it round cone} if $N$ is  convex, invariant with respect to $\Bbb R_{\ge 0}$-homotheties, Int$(N)\neq \emptyset$, and $N\cap N^{\circ} =\{ 0\}$. Here  $N^{\circ}:=-N$ is the opposite subset. Then $N^\vee := \{ \nu\in T^* : 
\langle \nu ,N \rangle \subset \Bbb R_{\le 0}\} $ is a round cone in the dual vector space $T^*$, and $(N^\vee )^\vee =N$. If $T$ is an $\Bbb R$-vector bundle over a topological space $X$, then  a {\it family of round cones}, or simply a {\it round cone}, in $T$ is a closed subspace $\CN\subset T$ such that each $\CN_x$, $x\in X$, is a round cone in $T_x$ which depends continuously on $x$. Then the family of dual cones $\CN^\vee$ is a round cone in $T^*$.  Notice that $\CN=\overline{\text{ Int}(\CN)}$, and Int$(\CN)_x =$ Int$(\CN_x )$. For  round cones $\CN,\CN' \subset T$  we write $\CN\Subset \CN'$ if $\CN\smallsetminus \{ 0\} \subset$ Int$(\CN')$. Sometimes we write $\CN_X$ instead of $\CN$; if  $Y\subset X$ is a subset (equipped with the induced topology), then  $\CN_Y \subset T|_Y$  is the pull-back of  $\CN_X$ to $Y$. 

Suppose now that our $X$ is a $C^1$-manifold; let $\pi_T : TX \to X$, $\pi_{T^*} : T^* X \to X$ be the tangent and the cotangent bundles.  We denote by $TX^\bullet$  the complement to the zero section of $TX$; for any $\CN\subset X$ set $\CN^\bullet := \CN\cap TX^\bullet$.

For the rest of the section we fix a round cone $\CN\subset TX$.

 A subset $P\subset X$ is said to be {\it $\CN$-invariant} if 
 for every $C^1$-arc $\gamma : [0,1] \to X$ with $\gamma (0)\in P$, $ 
(\frac{d}{ dt}\gamma )( [0,1])\subset \CN^\bullet$, one has $\gamma ([0,1]) \subset P$. 

\proclaim{\quad Lemma} (i)  $\CN$-invariantness is a local property:  $P\subset X$ is $\CN$-invariant if and only if for some (or any) open covering $\{ V_\alpha \}$ of $X$  each $P\cap V_\alpha \subset V_\alpha$ is $\CN_{V_\alpha}$-invariant.

(ii) If $P\subset X$ is $\CN$-invariant, then  $\bar{P}$, Int$(P)$ are $\CN$-invariant, and Int$(P)= Int(\bar{P})$. The subset Int$(P)$ is dense in $P$. More precisely, for any $C^1$-arc $\gamma : [0,1] \to X$ with $ (\frac{d}{ dt}\gamma )( [0,1])\subset \text{Int}(\CN)$ and $x\in \bar{P}$, one has $\gamma ((0,1])\subset $ Int$(P)$.

(iii) If $P,Q\subset X$ are $\CN$-invariant, then  $\overline{P\cap Q}= \bar{P}\cap \bar{Q}$.
\endproclaim

{\it Proof.} (i) is clear. (ii): When checking if a set is $\CN$-invariant, it suffices to consider only those $\gamma : [0,1] \to X$ which are embeddings. In this case,
 $\frac{d}{ dt}\gamma$ extends to a continuous vector field $\tau$ on a neighborhood of $\gamma ([0,1])$ which takes values in $\CN^\bullet$. Solving the differential equation, we get an open neighborhood $U\subset X$ of $x:=\gamma (0)$ and $C^1$-morphism $ U \times [0,1]\to X$, $(y,t)\mapsto g_t (y)$, such that $g_0 (y)=y$,  $\partial_t g_t (y)= \tau (g_t (y))$; then $g_t (x)=\gamma (t)$  and for each $t$ the map $g_t : U\to X$ is an open embedding. This  implies that $\bar{P}$, Int$(P)$ are $\CN$-invariant. 
 
Let $\CS$ be the set of all $C^1$-arcs $\gamma : [0,1] \to X$ such that $ (\frac{d}{ dt}\gamma )( [0,1])\subset \text{Int}(\CN)$. Then for each $n\ge 1$ the map $\CS \to X^{n+1}$, $\gamma \mapsto (\gamma (i/n))_{i=0,\ldots, n}$, has open image. 
 This observation shows that  for $x\in \bar{P}$ and $\gamma \in\CS$, $\gamma (0)=x$, one has $y:= \gamma (1)\in $ Int$(P)$. Indeed, to see that $y:= \gamma (1)\in $ Int$(P)$, notice that
$\gamma (1/2)\in\bar{P}$ by above;  thus for $z\in P$ sufficiently close to $\gamma (1/2)$   one can find $\gamma' \in\CS$ with $\gamma' (1/2)=z$ and $\gamma' (1)=y$. Then  the image of the set of all $\gamma'' \in\CS$ with $\gamma'' (1/2)=z$ by the map $\gamma'' \mapsto \gamma'' (1)$ is open; since  it lies in $P$, one has $y\in $ Int$(P)$, q.e.d.
 
It remains to check that Int$(P)=$ Int$(\bar{P})$, i.e., that every $y\in $ Int$(\bar{P})$ lies in Int$(P)$. Pick any $\gamma \in \CS$ such that $\gamma (1)=y$. Then  $\gamma (t)\in \bar{P}$ for any $t<1$ sufficiently close to 1, which implies that $\gamma (1)\in$ Int$(P)$, as we've seen.

(iii) It suffices to show that any $x\in \bar{P}\cap\bar{Q}$ lies in  $\overline{P\cap Q}$. 
Pick any $\gamma\in\CS$ with $\gamma (0)=x$; then $\gamma ((0,1])\subset P\cap Q$, and we are done.  \hfill$\square$.

 A subset $U\subset X$ is said to be {\it $\CN$-open} if it is open and $\CN$-invariant. Such subsets form
 a topology on $X$, called the {\it $\CN$-topology}.  Its closed sets, referred to as {\it $\CN$-closed} subsets, are the same as closed $\CN^{\circ}$-invariant subsets of $X$.

{\it Remark.} By (ii) of the above lemma, the map $U\mapsto \bar{U}$ is a 1--1 correspondence between the sets of $\CN$-open and $\CN^{\circ}$-closed subsets of $X$.

As the example of $X=S^1$ shows, the $\CN$-topology can be rather stupid. We will consider  instead  the  $\CN_U$-topologies for sufficiently small open $U$'s.

{\bf 2.3.} A $C^1$-function $f$ on an open subset $U\subset X$ is said to be 
{\it $\CN$-smooth} if $\tau (f)<0$ for every $\tau \in \CN_U^\bullet$, i.e., if $df(U)
\subset $ Int$(\CN^{\vee}_{U})$, or, equivalently, if there is $\CN'_U \Supset \CN_U$ such that each $U_{f<a}$ is $\CN'_U$-open.\footnote{If we call $\CN$  a light cone, then $f$ is a time function for it.}

\proclaim{\quad Lemma}  $\CN$-smooth functions exist locally.
\endproclaim

{\it Proof.} Take any $x\in X$. Choosing coordinates,  identify  a neighborhood $V$ of $x$ with that of $0\in\Bbb R^n$. Choose a round cone $N' \subset \Bbb R^n = T_x X$ such that $N'\Supset \CN_x$. Shrinking $V$, we can assume that $N'\Supset \CN_y$ for every $y\in V$. Then any non-zero linear function from ${N'}^\vee$ is $\CN$-smooth on $V$.  \hfill$\square$

A subset $P\subset X$ is said to be {\it $\CN$-small} if there is a neighborhood of $P$ which admits an $\CN$-smooth function. Of course, any subset of an $\CN$-small set is $\CN$-small.

\proclaim{\quad Lemma} Suppose that  $U\subset X$ is open and $\CN$-small. 

(i) Each $x\in U$ admits a base of open neighborhoods $\{ V_\alpha \}$ such that the $\CN_{V_\alpha}$-topology on $V_\alpha$ coincides with the topology induced by the $\CN_U$-topology on $U$. 

(ii) Suppose we have round cones $\CN''_U \Supset \CN'_U \supset \CN_U$ and an $\CN''_U$-open subset
$W\subset U$. Then for any $x\in U\cap \partial W$ the subsets $Q\smallsetminus W$, where $Q$ is an $\CN'_U$-neighborhood of $x$,  form a base of neighborhoods of $x$ in $U\smallsetminus W$.
\endproclaim

{\it Proof.} (i) Let $f$ be an $\CN$-smooth function on $U$. Choose a coordinate neighborhood $V$ of $ x$  and $N'$ as in the previous lemma so that the first coordinate $t_1$ equals $f$. For $\epsilon >0$ let $\tilde{\epsilon} := (\epsilon, 0,\ldots ,0)\in\Bbb R^n$. Set $R_\epsilon := $Int$(N') + \tilde{\epsilon}$, $V_\epsilon := (R_\epsilon )_{t_1 >-\epsilon}$; then $V_\epsilon$ form a base of neighborhoods of $0\in\Bbb R^n$. Let $\{ V_\alpha \}$ be the set of those $V_\epsilon $ that lie in $ V$; we leave it to the reader to check  the promised property.

(ii) Choose $V$ as above  and round cones $N', N'' \subset \Bbb R^n$ such that $\CN''_y \Supset N''\Supset N'\Supset \CN'_y$ for $y\in V$. The subsets $R_\epsilon \smallsetminus$ Int$(N'')$, where $R_\epsilon$ is as in (ii), form a base of neighborhoods of 0 in $\Bbb R^n \smallsetminus$ Int$(N'')$. Fix some $\epsilon'$ such that $V_{\epsilon'}\subset V$. Then $Q_\epsilon := (R_\epsilon \cap V_{\epsilon'})\cup W$ for $\epsilon \le\epsilon'$ are $\CN'_U$-open neighborhoods of $x$. For $\epsilon$  small enough,  $R_\epsilon \smallsetminus $ Int$(N'') \subset V_{\epsilon'}$, so
$Q_\epsilon \smallsetminus W = (R_\epsilon \smallsetminus $ Int$(N''))\smallsetminus W$; we are done.
\hfill$\square$

{\bf 2.4.} A subset $Z\subset X$ is called {\it $\CN$-special} if every $x\in \bar{Z}$ admits a neighborhood $V$ such that $Z\cap V$ is locally closed with respect to the $\CN'_V$-topology for some round cone $\CN'_V \Supset \CN_V$. If, in addition,
 $\bar{Z}$ is  compact and $\CN$-small, then $Z$ is called  {\it $\CN$-lens.}
 
\proclaim{\quad Lemma} (i) Suppose $Z$ is $\CN$-special. Then Int$(Z)= $ Int$(\bar{Z})$ is dense in $Z$. 

One recovers $Z$ from $U:=$ Int$(Z)$ as follows. For any $x \in \bar{U}\smallsetminus U$ there is its neighborhood $V$ such that for each arc $\gamma : [0,1]\to V$ with $\gamma (0)=x$, $(\frac{d}{ dt}\gamma )( [0,1])\subset \CN^{\circ\bullet}$, one either has $\gamma ((0,1])\subset U$ (then $x\in Z$), or $\gamma ((0,1])\subset V\smallsetminus \bar{U}$ (then $x\notin Z$).

(ii) If $Z$ is $\CN$-special, then it is $\CN\,\tilde{}$-special for every round cone $\CN\,\tilde{}$ which equals $\CN$ over Int$(Z)$. Same is true for the property of being $\CN$-lens.

(iii) The intersection of finitely many $\CN$-special subsets is $\CN$-special, same for $\CN$-lenses. Each point $x\in X$ admits a base of  neighborhoods formed by $\CN$-lenses.
\endproclaim

{\it Proof.} (i) Our assertions are local, so we can assume that $Z = W\cap Q$ where $W$ is $\CN$-open, $Q$ is $\CN$-closed. Then, by (ii) of the lemma in 2.2,  Int$(\bar{Z})\subset $ Int$(\bar{W})\cap$Int$(Q)=W\cap$Int$(Q)=$ Int$(Z)$, so Int$(Z)= $ Int$(\bar{Z})$. It is dense in $Z$ since Int$(Q)$ is dense in $Q$ by loc.~cit. The recipe for the recovery of $Z$ from Int$(Z)$ follows from loc.~cit.

  (ii) Since Int$(Z)$ is dense in $\bar{Z}$, the  cone $\CN'$ from the first paragraph of 2.4 works for $\CN\,\tilde{}$ as well.

  (iii) The first claim is evident. For the second one, choose $V$, $N'$  as in the proof of the first lemma in 2.3. For $a,b \in $ Int$(N' )$ set $Z^{N'}_{ab} := ($Int$(N' )-a)\cap (b-N') $.  Then  $Z^{N'}_{ab} $ form a base of neighborhoods of $0\in \Bbb R^n$. Those that lie in $ V$ are $\CN$-lenses, and we are done.  \hfill$\square$
  
For $U=$ Int$(Z)$ as in (i) above, we write $Z= U^+_{\CN}$, or simply $Z=U^+$. 

{\it Exercise.} Such an $U$ also equals Int$(Z^\circ )$ where $Z^\circ$ is a uniquely defined  $\CN^{\circ}$-special set. One has $Z\cap Z^\circ =U$. The map $Z\mapsto Z^\circ$ is a 1--1 correspondence between the sets of $\CN$- and $\CN^{\circ}$-special sets, which preserves the lenses and  extends the correspondence  from Remark in 2.2.

{\bf 2.5.} A complex of sheaves $F\in D(X)$ is said to be   {\it $\CN$-smooth} if  for every open $U\subset X$, a round cone $\CN'_U \Supset \CN_U$, and a $\CN'_U$-open $V\subset U$, one has $M_V (F)=0$, i.e., $i_{\bar{V}*}i_{\bar{V}}^* F \iso Rj_{V*} j_V^* F$ (see 2.1).  $\CN$-smooth $F$'s form a thick subcategory of $D(X)$.

{\it Remark.} One can consider an a priori weaker condition, demanding that $M_V (F)$ vanishes for all $V$ of type $X_{f<a}$, where $f$ is an $\CN$-smooth function. We will see that this property is equivalent to $\CN$-smoothness.

\proclaim{\quad Lemma} If $F\in D(X)$ is $\CN$-smooth and  $Z\subset X$ is $\CN$-special, then the complexes $Ri_{Z*} Ri_Z^! F$ and $i_{Z^\circ !} i^*_{Z^\circ}F $ are $\CN$-smooth, and the canonical morphism $i_{Z^\circ !} i^*_{Z^\circ}F \to Ri_{Z*} Ri_Z^! F$ (which equals $id_F$ on Int$(Z)=Z\cap Z^\circ$) is a quasi-isomorphism. \endproclaim

{\it Proof.} Our assertions are local, so we can assume that $Z=U\smallsetminus U'$, where $U$, $U'$ are $\CN'$-open for some $\CN'\Supset \CN$. Then $Z^\circ = \bar{U}\smallsetminus \bar{U}'$.

  Let us check first that $G:= Rj_{U*} j_U^* F  $ is $\CN$-smooth. For any $\CN'$-open $V$ one has $Rj_{V*} j_V^* G = Rj_{U\cap V*} j_{U\cap V}^* F$, which is     $i_{\overline{U\cap V}*}i_{\overline{U\cap V}}^* F$  by $\CN$-smoothness of $F$.  By the same reason, one has $i_{\bar{V}*}i_{\bar{V}}^*G=i_{\bar{V}*}i_{\bar{V}}^*i_{\bar{U}*}i_{\bar{U}}^* F= i_{\bar{U}\cap\bar{V}*}i_{\bar{U}\cap\bar{V}}^*F$. Since $\overline{U\cap V}=\bar{U}\cap\bar{V}$ by 2.2(iii), one has $i_{\bar{V}*}i_{\bar{V}}^*G\iso Rj_{V*} j_V^* G$; we are done. Also, since $F$ is $\CN'$-smooth, one has  $i_{U^\circ !} i^*_{U^\circ} F= i_{\bar{U}*} i^*_{\bar{U}}=G$.

Now one has $Ri_{Z*} Ri_Z^! F = Rj_{U*} j_U^* \CC one (F\to Rj_{U'*} j_{U'}^* F)[-1]$ and  
 $i_{Z^\circ !} i^*_{Z^\circ}F = i_{\bar{U}*}i_{\bar{U}}^* \CC one (F\to i_{\bar{U}'*}
   i_{\bar{U}'}^*F)[-1]$. Therefore, by the above,  the complexes
     $Ri_{Z*} Ri_Z^! F$ and  $i_{Z^\circ !} i^*_{Z^\circ}F$  are equal and $\CN$-smooth, q.e.d.   \hfill$\square$

{\it Exercise.} If $Z$ is $\CN$-special, then any
 $\CN$-smooth complex  supported on $\bar{Z}$ is determined uniquely (up to a unique isomorphism) by its restriction to Int$(Z)$.

\proclaim{\quad Proposition} For $F\in D(X)$ the next properties are equivalent:

(i) $F$ is $\CN$-smooth;

(ii) $F$ satisfies the property from Remark above;

(iii) For each $\CN$-lens $Z$ one has $R\Gamma  (Z, Ri^!_Z F)=0$;

(iv) For each $\CN$-lens $Z$ one has $R\Gamma_c  (Z^\circ , i^*_{Z^\circ} F)=0$.
\endproclaim

{\it Proof.}  (iii)$\Rightarrow$(i): Pick any $U$, $V$, $\CN'_U$ as in the beginning of 2.5; let us show that $M_V (F) =0$. Our problem is local, so we can assume that $\bar{U}$ is compact and $\CN$-small.
It suffices to check that $R\Gamma (U', M_V (F))=0$ for every open $U' $ such that $\bar{U}' \subset U$.
Since $M_V (F) = i_{\partial V *} i_{\partial V}^* R i^!_{U\smallsetminus V} F$, one has $R\Gamma (U', M_V (F))=R\Gamma (U' \cap\partial V ,  i_{\partial V}^* R i^!_{U\smallsetminus V} F )\buildrel\sim\over\leftarrow \limright R\Gamma (Z, Ri^!_Z F)$, where $Z$ runs the directed set $\CT$ of open neighborhoods of $U'\cap\partial V$ in $U'\smallsetminus V$. Choose $\CN''_U$ such that $\CN'_U \Supset \CN''_U \Supset \CN_U$. Then, by (ii) of the second lemma in 2.3,\footnote{Sorry for  the discrepancy of  notation: the present $\CN'_U$, $\CN''_U$, $V$ are $\CN''_U$, $\CN'_U$, $W$ in loc.~cit.}  those $Z$ which can be represented as $Q\smallsetminus V$, where $Q\subset U$ is an $\CN''_U$-open subset, form a cofinal subset in $\CT$. Such a $Z$ is an $\CN_U$-lens, hence $R\Gamma (Z, Ri^!_Z F)$ vanishes; we are done.

(iv)$\Rightarrow$(i): We use the dual argument, handling accurately the inverse limits: 

Let $U$, $V$, $\CN'_U$ be as above. To prove that $M_V (F)$ vanishes, it suffices to show that for every $U'$ as above the complex $K_{U'} := R\Gamma_c (U', M_V (F))$ is acyclic. One  can assume that $F$ is a complex of flabby sheaves; then  $K_{U'}=\Gamma_c (U' \cap\partial V ,  i_{\partial V}^!   i^*_{\bar{V}} F ) $ is the subcomplex of $\Gamma_c (U' \cap \bar{V},  i^*_{\bar{V}} F ) $ formed by sections that vanish on $U'\smallsetminus \partial V$. Since $K_{U'}= \cup K_{U''}$, the union is taken with respect to all open $U''$ with $\bar{U}'' \subset U'$, it suffices to show that each map $H^0 K_{U''}\to H^0 K_{U'}$ vanishes.
 
 Let $\CT'$ be the collection of all subsets $Z'\subset U''$ which contain $U'' \cap\partial V$ and can be represented as $P\cap \bar{V}$, where $P\subset U$ is an ${\CN''_U}^{\circ}$-open subset. 
  By (ii) of the second lemma in 2.3,\footnote{The present $\CN'_U$, $\CN''_U$, $U\smallsetminus \bar{V}$ are ${\CN''_U}^{\circ}$, ${\CN'_U}^{\circ}$, $W$ of loc.~cit.} such $Z'$'s are cofinal among the open neighborhoods of $U'' \cap\partial V$ in $U'' \cap\bar{V}$. So there is a decreasing sequence $Z'_0 \supset Z'_1 \supset \ldots$ in $\CT'$ with $\cap Z'_i = U''\cap\partial V$. 
    Set $L_i :=\Gamma_c (Z'_i , i^*_{Z'_i} F)$; then $L_0 \supset L_1 \supset \ldots$, and $K_{U''}=\cap L_i$. Since $L_i = R\Gamma_c (Z'_i , i^*_{Z'_i}F )$ and $Z'_i$ is an $\CN^{\circ}$-lens, each $L_i$  is acyclic. So $\hat{L}:= \limleft L_0 /L_i$ is also acyclic. Notice that $\hat{L}$ is the subcomplex of $\Gamma (Z'_0 \smallsetminus \partial V,F)$ that consists of those sections $\gamma$ that the closure of supp$(\gamma )$ is contained in $Z'_0$. Since $F^i$ are flabby, any such $\gamma$ can be extended to a  section $\tilde{\gamma}\in \Gamma (U' \cap \bar{V}, i^*_{\bar{V}}F)$ whose  support is compact and lies in supp$(\gamma )\cup \partial V$
  
Now let $\alpha \in K_{U''}^0$ be any cycle. Then $\alpha = d\beta$ for some $\beta \in L^{-1}_0$. The image $\hat{\beta}$ of $\beta$ in $\hat{L}^{-1}$ is a cycle, hence $\hat{\beta}=d\gamma$ for some $\gamma \in \hat{L}^{-3}$. Pick $\tilde{\gamma}$ as above. Then $\beta' :=\beta - d\tilde{\gamma} \in K_{U'}^{-1}$ and $d\beta' =\alpha$, i.e., $\alpha$ is exact in $K_{U'}$, q.e.d.

(i)$\Rightarrow$(iii),(iv): By Exercise, for  $\CN$-smooth $F$  assertions (iii) and (iv) are equivalent; we check (iii).
By Lemma, it suffices to show that for any $\CN$-smooth $F$ whose support $S$ is compact and $\CN$-small, one has $R\Gamma (X,F)=0$. Since $S$ is compact,  one can choose a $\CN$-smooth function $f$ on a neighborhood of $S$ such that $f(S)\subset (0,1)$. Now the lemma from 2.2 (for $X$ in loc.~cit.~equal to $S$, $V_t := S_{f<t}$, and the sheaf $i^*_S F$) does the job.

 (ii)$\Rightarrow$(i): Let $U$, $V$, $\CN'_U$ be as in the beginning of 2.5.  Pick $x\in\partial V$; let us check that $M_V (F)_x =0$. The problem is local, so, as in the proofs in 2.3, we can assume that $U$ is an open subset of $\Bbb R^n$, $x=0$, and there is a round cone $N' \subset \Bbb R^n$ such that $N'_U \Supset \CN_U$ and $V$ is $N'_U$-invariant. 
 
Pick any $\ell \in $ Int$({N'}^\vee )$; for $y\in \Bbb R^n$, $a\in\Bbb R$ set $Z_{ay}:=  \ell_{\ge a}\cap$Int$(N'+y)\subset \Bbb R^n$. One easily finds a $C^\infty$-function $\phi (v,t)=\phi_t (v)$ on $\Bbb R^n \times [0,1)$ such that $\phi_0 =\ell -a$, each $\phi_t$ is $N'_{\Bbb R^n}$-smooth (i.e., $d\phi_t $ takes values in Int$(N^\vee )$) and decreases as a function of $t$, and $\cup_t \Bbb R^n_{\phi_t <0} =$ Int$(N' +y)\cup \Bbb R^n_{\ell <0}$. Applying the lemma from 2.2 to $V_t =\Bbb R^n_{\phi_t <0}$ (and $V_1 = $ Int$(N' +y)\cup \Bbb R^n_{\ell <a}$), we see that 
$R\Gamma (Z_{ay}, Ri_{Z_{ay}}^! F)=0$ whenever $\bar{Z}_{ay}\in U$.

Now $M_V (F)_x = \limright \CC one (R\Gamma  (Z_{ay}, Ri_{Z_{ay}}^! F)\to 
 R\Gamma  (V\cap Z_{ay}, Ri_{Z_{ay}}^! F))[-1]$, where $y\in -$Int$(N')$ tends to 0 and $a<0$ is  such that  $\bar{Z}_{ay}\subset U$ for some  $y\in -$Int$(N')$ (cf.~the proof of  (iii)$\Rightarrow$(i)). Since $V\cap Z_{ay}$ admits an (open) hypercovering by  $Z_{av}$'s,   $v \in V\cap Z_{ay}$, the vanishing of $R\Gamma (Z_{ay}, Ri_{Z_{ay}}^! F)$  implies that $M_V (F)_x =0$, q.e.d.   \hfill$\square$

As in \cite{KS} 5.1, one defines the {\it micro-support}  of $F\in D(X)$ as the smallest closed subset $SS(F)\subset T^* X$ such that for every open $U\subset X$ and a $C^1$-function $f$ on $U$ with property $df (U)\cap SS(F)=\emptyset$, the Morse complexes $M_{f=a} (F)\in D(U)$ vanish. So condition (ii) of the proposition can be restated as $SS(F)\cap$Int$( \CN^\vee )=\emptyset$.

{\bf 2.6.} {\it Notation.}  A {\it partition} $\{ Z_\alpha \}$ of a topological space $Z$ is always 
assumed to be finite.  We say that a partition $\{ Z_\alpha \}$  is {\it locally closed} if each 
$Z_\alpha$ is locally closed. 
A {\it filtration} on $Z$ always means finite filtration 
 $\emptyset = Z_{-1} \subset \ldots \subset Z_m = Z$; such a filtration defines a partition 
 $\{ Z_i \smallsetminus Z_{i-1}\}$. Our filtration is {\it closed} if $Z_i$'s are closed subsets. 
 We say that a  partition $\{ Z_\alpha \}$ 
 {\it comes from a closed filtration} if there is a  closed filtration such that the corresponding 
 partition coincides with $\{ Z_\alpha \}$. Such a filtration amounts to
  a linear ordering on the set of 
 indices $\alpha$ such that for every $\alpha$ the subset $\bar{Z}_\alpha \smallsetminus Z_{\alpha}
  $ is closed and lies in the union 
 of subsets $Z_{\alpha'}$, $\alpha' < \alpha$. A locally closed partition $\{ Z_\alpha \}$ is a 
 {\it stratification} if for every $\alpha$ the closure $\bar{Z}_\alpha$ is a union of parts 
 $Z_{\alpha'}$ (the parts $Z_\alpha$ are called then  {\it strata}). 
The set Par$(Z)$ of  partitions is ordered with respect to the
refinement: one has $\{ Z_\alpha \} \le \{ Z_\beta \}$ if for every $\beta$ there is a 
(necessary unique) $\alpha =\alpha (\beta )$ such that 
$Z_\beta \subset Z_{\alpha (\beta )}$.
 Every finite subset of Par$(Z)$ admits the least upper bound, so Par$(Z)$ is directed. 
The class Par$^{c}(Z)$ of partitions that come from a closed
 filtration satisfies the transitivity property: 
if $\{ Z_{\alpha}\} \le \{ Z_{\beta}\} $ are such that $\{ Z_{\alpha }\}\in$ Par$^{c}(Z)$, 
 and
for each $Z_{\alpha}$ the partition $\{ Z_{\beta_{\alpha}}\}$ of $Z_{\alpha}$
 formed by those $Z_{\beta}$ that lie in
$Z_{\alpha}$ belongs to Par$^{c}(Z_{\alpha})$, then $\{Z_{\beta}\}\in$ Par$^{c}(Z)$.\footnote{To see this, consider the lexicographic ordering of $\beta$'s.}
 This implies that  Par$^{c}(Z)$
  is closed with respect to the least upper bound of
finite subsets, so 
it is a directed subset of Par$(Z)$. Same is true for the larger class of locally closed partitions.

{\bf 2.7.} We return to our situation. Let $Z$ be an $\CN$-lens. A   partition $\{ Z_\alpha \}$ of $Z$ is said
 to be {\it $\CN$-special} if it comes from a closed filtration and each $Z_\alpha$ is 
 $\CN$-special (i.e., is an $\CN$-lens). The filtration is referred to as {\it $\CN$-special}
 filtration. The  set   Par$^{sp}_{\CN}(Z)\subset $ Par$^{c}(Z)$ 
 of $\CN$-special partitions is closed with respect 
 to the least upper bound of finite subsets, hence is directed, and satisfies the above transitivity property.
 
 \proclaim{\quad Lemma} (i) Suppose we have an open covering of $\bar{Z}$ in $X$. Then 
 one can find $\{ Z_\alpha \}\in$ Par$^{sp}_{\CN }(Z)$ such that each $Z_\alpha$  lies in 
 some open subset of the covering.  
 
 (ii) If, in addition, we have a finite collection of $\CN$-lenses $ Z'_i $, then there is 
 $\{ Z_{\alpha }\} $  which is 
 finer than all the partitions $\{ Z\cap Z'_{i}, Z\smallsetminus Z'_{i} \}$.
 \endproclaim
 {\it Proof.} (i) Replacing $X$ by a neighborhood of $\bar{Z}$, we can assume that there
  is a global $\CN$-smooth function $f$. Refining the open covering $\{ U_\beta \}$, we can assume that
   for each  $U_\beta$  the intersection $Z\cap U_\beta$ is locally closed in $U_\beta$ with respect to the 
   $\CN_{U_\beta}$-topology. 
 
 Take any $x\in \bar{Z}$, and choose $U_\beta$ that contains $x$.   Then one can find a neighborhood 
 $P_x$ of $x$ such that $\bar{P}_x \subset U_\beta$, $P_{x}$  is locally closed in
  $\CN_{U_\beta}$-topology, and $f(x)< f(y) $ for every $y\in \bar{P}_x \smallsetminus P_x$. 
  To see this, we proceed as in the proof of part (iii) of the lemma in 2.4, assuming in addition that 
   $df (x)\in {N'}^\vee$. Then $df$ takes values in ${N'}^\vee$ on a neighborhood of  $x$, and one can choose $P_x$ among $Z^{N'}_{ab}$'s from loc.~cit.~(with $a,b$ small, and $b$ much smaller than $a$).
  
 For every $a\in\Bbb R$ the open subsets Int$(P_x )$, $x\in \bar{Z}_{f=a}$, cover the compact set $\bar{Z}_{f=a}$. Choose a finite subset $P^a_1 ,\ldots , P^a_{N(a)}$ of these $P_x$'s such that Int$(P^a_i )$ cover $\bar{Z}_{f=a}$. There is an open interval $ I_a \ni a$ such that $f^{-1} (I_a )\subset \mathop\cup\limits_i P^a_i$, $f (\mathop\cup\limits_i (\bar{P}^a_i\smallsetminus P^a_i ))\cap I_a = \emptyset$.
 
 The intervals $I_a$ cover  $f(\bar{Z})$. Let $I_{a_1 },\ldots, I_{a_m}$ be a finite subcovering; we
  order $a_j$'s so that the top points of $I_{a_j}$'s decrease. Consider the collection of the subsets $P_{ij}:= P^{a_j}_i$. 
   Ordering the indices lexicografically, denote these subsets as $P_1, \ldots , P_N$. Then
   $(\bar{P}_{i}\smallsetminus P_{i})\cap Z\subset \mathop\cup\limits_{k<i}P_{k}$, and
  $Z_i := \mathop\cup\limits_{k\le i} (P_k \cap Z)$ is an $\CN$-special filtration of the kind
  we look for.

    (ii)  Since Par$^{sp}_{\CN}(Z)$ is directed, it suffices to consider the case of a single $Z'$. We modify 
    the above procedure as follows. 
    
    Refining $\{ U_{\beta}\}$, we can assume that each  $Z'\cap U_\beta$ is locally closed in $U_\beta$ with respect to the 
   $\CN_{U_\beta}$-topology, so $Z'\cap U_{\beta} = Q_{\beta}^{1}\smallsetminus Q_{\beta}^{2}$ 
   where $Q_{\beta}^{i}\subset U_{\beta}$ are $\CN_{U_{\beta}}$-closed subsets. 
   For $x\in \bar{Z}$ we  pick  $P_{x}\subset U_{\beta}$ as above, and 
    set $P_{x 2} := P_{x}\smallsetminus Q_{\beta}^{1}$, 
    $P_{x1}
    := P_{x}\cap (Q_{\beta}^{1}\smallsetminus Q_{\beta}^{2})$,  $P_{x0}:=P_{x}\cap Q_{\beta}^{2}$.
    
   Choose from all $P_{x}$'s the subsets  $P_{ij}$ as above; for $\ell=0,1,2$ let $P'_{i\ell j}\subset P_{ij}$
   be the corresponding $P_{x\ell}$. Ordering the indices lexicographically, denote the subsets $P'_{i\ell j}$ by 
   $P'_{1},\ldots , P'_{M}$.  Then
   $(\bar{P}'_{i}\smallsetminus P'_{i})\cap Z\subset \mathop\cup\limits_{k<i}P'_{k}$, and
  $Z_i := \mathop\cup\limits_{k\le i} (P'_k \cap Z)$ is an $\CN$-special filtration of the kind
  we look for. 
     \hfill$\square$
 
{\bf 2.8.} The class of {\it $\CN$-constructible} subsets of $X$ comes from that of $\CN$-lenses by closing
it with respect to finite unions, intersections and differences. 

{\it Remark.} $\CN$-constructibility is essentially a local property. Namely, a subset $Y\subset
X$ is $\CN$-constructible if (and only if) $\bar{Y}$ is compact and for some  (or every) open 
covering $\{ U_{i}\}$ of $X$
each $Y\cap U_{i}$ belongs to the Boolean algebra of the subsets of $U_{i}$ generated by the $\CN_{U_{i}}
$-special subsets.\footnote{To see this, pick for any $x\in X$ a neighborhood $V_{x}$ 
which is an $\CN$-lens such that
$\bar{V}_{x}$ lies in some $U_{i}$, and choose finitely many of $V_{x}$'s that cover $\bar{Y}$.
 Then each $Y\cap V_{x}$ is $\CN$-constructible, and their union equals $V$, so $V$ is $\CN$-constructible.}

Let $Y$ be an $\CN$-constructible set.
Denote by Par$^{\prec}_{\CN}(Y)$ the set of partitions $\{ Y_{\alpha}\}$ of $Y$ such that 
each $Y_{\alpha}$ is an $\CN$-lens.  It carries a partial 
order $\prec$ which is finer than the order $<$ coming from Par$(Y)$: namely,
  $\{ Y_{\alpha} \} \preceq \{ Y_{\beta}\}$ if for each $\alpha$ those $Y_{\beta}
$ that lie in $Y_{\alpha}$ form an $\CN$-special partition of $Y_{\alpha}$ (see 2.7).

\proclaim{\quad Proposition} Par$^{\prec}_{\CN} (Y)$ is a non-empty directed poset. 
\endproclaim

{\it Proof.} $Y$ admits a partition $\{ Y_{\beta} \}$  such that each $Y_{\beta}$ can be 
written as $Z\smallsetminus \cup Z'_{i}$, 
where $Z$, $\{ Z'_{i} \}$ are $\CN$-lenses. By 2.7(ii), $Y_{\beta}$ admits a 
partition by $\CN$-lenses $\{ Y_{\alpha\beta}\}$. Then all $Y_{\alpha\beta}$  together
form a partition of $Y$ by $\CN$-lenses, so Par$^{\prec}_{\CN} (Y)\neq\emptyset$.

It remains to check that Par$^{\prec}_{\CN}(Y)$ is directed. Suppose we have
 $\{ Y_{\alpha }\},\{ Y_{\beta}\}\in$ Par$^{\prec}_{\CN}(Y)$; let us find their upper bound
  with respect to $\prec$. Let $\{ Y_{\gamma}\}$ be their least upper bound for $<$.
  Then $\{ Y_{\gamma} \} \in $ Par$^{\prec}_{\CN }(Y)$. Applying 2.7(ii) to each $Z=
  Y_{\alpha}$ and $Z'_{i}=Y_{\gamma}$, we find $\{ Y_{\delta} \}\in$ Par$^{\prec}_{\CN}(Y)$ 
   such that $\{ Y_{\alpha}\} 
  \preceq \{ Y_{\delta }\}$ and $\{ Y_{\gamma}\} \le \{ Y_{\delta}\}$. Applying 2.7(ii) to each 
  $Z=Y_{\beta}$ and $Z'_{i }=Y_{\delta}$, we find $\{ Y_{\epsilon}\} \in$ 
  Par$^{\prec}_{\CN }(Y)$
  such that $\{ Y_{\beta}\} \preceq \{ Y_{\epsilon}\}$ and $\{ Y_{\delta}\}\le \{ Y_{\epsilon}\}$.
  The latter inequalities imply that $\{ Y_{\delta}\} \preceq\{ Y_{\epsilon}\}$, 
  so $\{ Y_{\epsilon}\}$ is an
  upper bound of $\{ Y_{\alpha}\}, \{ Y_{\beta}\}$ for $\prec$, q.e.d.
  \hfill$\square$
 
 \bigskip

\centerline {\bf 3. $ K$-spectra of constructible sheaves and the  $\varepsilon$-factorization.}

\medskip

In this section we explain the constructions mentioned in 0.1(a).

{\bf 3.1.} For the next notions and facts see  \cite{KS} 8.2 and references therein.  

From now on our $X$ is a real analytic manifold. Recall that a subset $Z$ of $X$
is  {\it subanalytic} if every point of $X$ admits an open neighborhood
$U$ such that $Z\cap U$ belongs to  the Boolean algebra of subsets of $U$  generated by the images
of proper real analytic maps $Y\to U$,  $Y$ are real analytic manifolds. If $X$ is compact, then one can take for $U$ the whole $X$.
The property of being subanalytic is local; the closure of a subanalytic subset is subanalytic; 
every closed subanalytic subset can be represented as the image of a proper map as above. 

We use the terminology from 2.6. A  partition $\{ Z_\alpha \}$  of a subset of $X$ is  {\it subanalytic} if each $Z_\alpha$ is subanalytic
   in $X$. 
   Notice that if $\{ Z_\alpha \}$ is a subanalytic stratification of a locally closed subset $Z\subset X$ such that each $Z_\alpha$ is 
   smooth equidimensional, then $\{ Z_\alpha \}$ comes from a closed filtration. If $\bar{Z}$ is compact, then
 subanalytic  stratifications 
$\{ Z_\alpha \}$ with each $Z_\alpha$ smooth, equidimensional, and contractible, are cofinal in 
the set Par$^a (Z)$ of subanalytic partitions. Moreover, those that are homeomorphic to a simplicial 
decomposition still do.

Let $R$ be an associative unital ring\footnote{In fact, $R$ can be any associative DG algebra, or even a small 
DG category, as in \cite{Dr}.}
 and $\CP\subset D(R)$ be a full triangulated subcategory in the derived category of complexes of $R$-modules. 
 
For a locally closed subanalytic subset $Z$ a  {\it $\CP$-constructible complex} on $Z$ 
 is a complex of sheaves of $R$-modules  $F$  on $Z$ with the following
 properties: (a) every point of $\bar{Z}$ admits a neighborhood $U$ such that for
 some  $\{ Z_\alpha \} \in $ Par$^a (Z\cap U)$  the restrictions of the cohomology sheaves 
 $H^i (F)|_{Z_\alpha}$ are all constant, and (b) for each $x\in Z$ one has $F_x \in \CP$.  If $\bar{Z}$ is compact, then
 one can take $U=X$.
 Such complexes form a full triangulated subcategory $\CP(Z)$ of the derived category of complexes of 
 sheaves of $R$-modules on $Z$. We assume that $\CP$ is essentially small;  then such is $\CP(Z)$. 
 We have  the K-theory spectra $K :=K(\CP)$ and $K (Z):=K(\CP(Z))$ (see 1.8).
 
If   $f: X\to Y$ a proper real analytic map, then $Rf_{*}$ sends $\CP (X)$ to $\CP(Y)$, so we have the homotopy
morphism of spectra $Rf_{* }: K(X)\to K(Y)$. If  
 $T\subset Y$ is a locally closed subanalytic subset that contains $f(Z)$, and $g:= f|_Z : Z\to T$, then 
 the functors $Rg_!$, $Rg_*$ send $\CP(Z)$ to $\CP(T)$, so  we have  homotopy morphisms of spectra $Rg_! ,Rg_* : K (Z)\to K (T)$. 
 If $Y$ is a point, then we write $R\Gamma =R\Gamma (Z,\cdot )$ instead of  $Rf_* $, and $R\Gamma_c =R\Gamma_c (Z,\cdot )$ for $Rg_!$.

{\it Principal example:} $\CP=D^{per\! f}(R)$ is formed by perfect $R$-complexes, so $K =K(R)$, and $\CP(Z)$ consists of perfect constructible $R$-complexes of sheaves.

\medskip

{\bf 3.2.} Let $Z\subset X$ be a locally closed  subanalytic subset with $\bar{Z}$ compact,
 $\{ Z_\alpha \}$ be its  locally closed subanalytic partition.
 Denote by $i_\alpha $ the embeddings $ Z_\alpha \hra Z$. Set $K(\{ Z_{\alpha}\}):= \mathop\prod\limits_\alpha K (Z_\alpha  )
  = K (\mathop\sqcup\limits_\alpha Z_\alpha )$.

\proclaim{\quad  Lemma}  The morphisms $$
K (Z )
\overset{(i_\alpha^*)}\to{\underset{\Sigma i_{\alpha !}}
\to\rightleftarrows} K(\{ Z_{\alpha} \})\tag 3.2.1$$ are homotopically mutually inverse weak equivalences of spectra. 
The arrows in (3.2.1) are transitive, in the evident sense, with respect to the refinement of partitions.

\endproclaim

{\it Proof.}  (a) The transitivity statement  is clear. Since $(i_\alpha^*)$ is the  left inverse to $\Sigma i_{\alpha !}$, it suffices to check that either of the arrows in (3.2.1) is a weak equivalence. 

(b) It suffices to consider the case when our partition comes from a closed filtration. Indeed, for an arbitrary $\{ Z_\alpha \}$ one can find  $\{ Z_\beta \}$ that comes from a closed filtration, which is finer than $\{ Z_\alpha \}$. Then for any $\alpha $ the subset $\{ Z_{\beta_\alpha }\}$ of those $Z_\beta$ that lie in $Z_\alpha$ is a partition of $Z_\alpha$ that comes from a closed filtration.  By the transitivity,  our  assertion for $(Z, \{ Z_\alpha \} )$ follows from those for $(Z,\{ Z_\beta \} )$ and $(Z_\alpha ,\{ Z_{\beta_\alpha }\})$.

(c) Assuming that our partition comes from a closed filtration, let us show
that the composition $(\Sigma i_{\alpha !})(i_\alpha^*)$ is homotopic to the identity. 
A  closed filtration $Z_n$ on $Z$  yields a natural filtration on every sheaf $F$ on $Z$
with $\gr F =\oplus i_{\alpha !}i_\alpha^* F$. Thus  the identity functor of $\CP (Z)$ carries a canonical filtration with $\gr $ equal to $ \oplus  i_{\alpha !}i_\alpha^*$. By 1.8(iv), this implies that  $
(\Sigma i_{\alpha !})(i_\alpha^*)$ is naturally homotopic to $id_{K (Z )}$, q.e.d.   \hfill$\square$ 

{\it Remark.} Since the composition $K(\{ Z_\alpha \}) \buildrel{\Sigma\, i_{\alpha !}}\over\lra K(Z)\buildrel{R\Gamma_c}\over\lra K$ has components $R\Gamma_c : K(Z_\alpha )\to K$, we see that
$R\Gamma_c = \Sigma_\alpha \, R\Gamma_c i^*_\alpha : K(Z)\to K$.

\medskip

 {\bf 3.3.}  This subsection will not be used in the sequel; the reader can skip it.

 For $Z$ as in 3.2
  and a  spectrum $C$ let $\CC ons (Z,C)$ be the spectrum of constructible $C$-valued functions on $Z$. Precisely, for each  $\{ Z_\alpha \} \in$ Par$^a (Z)$ we have a spectrum $C^{\{Z_\alpha \}}$ of $C$-valued functions on $Z$ constant along $ Z_\alpha $'s;  if $\{ Z_\beta \}$ is a finer partition, then the projection $ \{ Z_\beta \} \twoheadrightarrow \{ Z_\alpha \}$ yields an  embedding $C^{\{Z_\alpha \}}\subset C^{\{Z_\beta \}}$. One has $\CC ons (Z,C):=\mathop\cup\limits_{\text{Par}^a (Z)} C^{\{Z_\alpha \}}$; thus $\pi_i \CC ons (Z,C)$
  is the group of $\pi_i C$-valued constructible functions on $Z$.

\proclaim{\quad Lemma} There is a canonical weak equivalence of spectra 
$$K (Z)\iso \CC ons (Z, K ) .\tag 3.3.1$$
\endproclaim

{\it Proof.} Let $S\subset $ Par$^a (Z)$ be the subset of those $\{ Z_\alpha \}$ that come from a closed
filtration
 and each $Z_\alpha$ is contractible. It is cofinal in Par$^a (Z)$, hence
  $\CC ons (Z, K )=\mathop\cup\limits_{S} K^{\{ Z_\alpha \}} $  (see 3.1). For $\{ Z_\alpha \} \in S$ let $\CP (Z,\{ Z_\alpha \})\subset \CP (Z)$ be the full DG subcategory of complexes $F$ such that all sheaves $H^a (F)|_{Z_\alpha}$ are constant; set $K (Z,{\{ Z_\alpha \}}):= K(\CP (Z,{\{ Z_\alpha \}}))$. Since $\CP (Z) 
=\mathop\cup\limits_{S }\CP (Z,\{ Z_\alpha \})$, one has 
$K (Z )=  \mathop\cup\limits_{S}K (Z ,\{ Z_\alpha \}) $. Since $S$ is a directed  poset,  one can replace $\cup$ by the homotopy colimit.

Choose a natural flabby resolution $F\to \tilde{F}$. For $\{ Z_\alpha \} \in S$  consider  DG functors 
$\CP (Z,\{ Z_\alpha \})\overset{\psi_{\{ Z_\alpha\}}}\to{\underset{\chi_{\{ Z_\alpha\}}}
\to\rightleftarrows} \CP^{\{ Z_\alpha \}}$, 
where $\psi_{\{ Z_\alpha\}}(F):= (\Gamma (Z_\alpha , i^*_\alpha \tilde{F} ))$, and 
$\chi_{\{ Z_\alpha\}} ((L_\alpha )):=\oplus  i_{\alpha !} 
L_{\alpha Z_\alpha}$.\footnote{Here $ L_{\alpha Z_\alpha}$ is the constant sheaf on $Z_\alpha$ with fibers
 $L_\alpha $.} By the argument from the proof in 3.2, the 
 morphisms of the K-spectra $K (Z , \{ Z_\alpha \})\rightleftarrows  K^{\{ Z_\alpha  \} }$ are homotopically mutually inverse weak equivalences.

If $\{ Z_\beta \} \in S$ is finer than $\{ Z_\alpha \}$, then
 for $F\in \CP (Z,\{ Z_\alpha \})$ one has an evident natural quasi-isomorphism 
 $\psi_{\{ Z_\alpha \}} (F)\to \psi_{\{ Z_\beta \}} (F)$ in $\CP^{\{ Z_\beta \}}$.
  It satisfies the transitivity property, so we get a weak equivalence
   $\psi : $ hocolim$_{S} \, K (Z,\{ Z_\alpha \})\to \CC ons(Z,K )$, which is (3.3.1). \hfill$\square$

\medskip
{\bf 3.4.}  We need a version of 3.2 in the presence of a round cone $\CN\subset TX$.

We use a subanalytic version of some structures from \S 2. Namely, an $\CN$-lens $Z$ (see  2.4) is said to be  subanalytic if  one can choose $V$ and $\CN'_{V}$ as in loc.~cit.~so that 
 $Z\cap V=U_{1} \smallsetminus U_{2}$ where $U_{i}$ are $\CN'_{V}$-open subanalytic subsets.
  As in 2.8, subanalytic $\CN$-lenses give rise to subanalytic $\CN$-constructible subsets. The results
  of 2.4, 2.7, 2.8 remain valid (together with the proofs) in the subanalytic setting. {\it From now on, all $\CN$-lenses
  and $\CN$-constructible sets are tacitly assumed to be subanalytic.}
  
 The complexes  $F$ with $SS(F)\cap$Int$(\CN^{\vee} )=\emptyset$
 form a thick subcategory $\CP (X)^{\CN} \subset\CP (X)$;
 set $K^\CN (X) := K(\CP (X)^{\CN})$. For an $\CN^{\circ}$-lens $Z$ let $\CP (Z)^{\CN}$ be a thick subcategory of $\CP (Z)$
  formed by those complexes $F$ that $i_{Z!}F\in\CP (X)^{\CN}$; set $K^\CN (Z):=K(\CP (Z)^{\CN})$.
  By 2.5, the functor $i_{Z}^{*}$ sends $\CP (X)^{\CN}$ to $\CP (Z)^{\CN}$; similarly, if $i: Z_{1}\hra 
  Z_{2}$ are
     $\CN^{\circ}$-lenses, then $i_{!}$ and
  $i^{*}$ interchange $\CP (Z_{i})^{\CN}$. Therefore we have the  morphisms 
  of spectra $K^\CN (X) \overset{i_Z^*}\to{\underset{ i_{Z !}}\to\rightleftarrows} 
   K^\CN (Z)$, $K^\CN (Z_{2})
\overset{i^*}\to{\underset{ i_{ !}}
\to\rightleftarrows}  K^\CN (Z_{1})$.

\proclaim{\quad Lemma}  Let $Z$ be an $\CN^{\circ}$-lens, 
$\{ Z_{\alpha}\}$ its partition by $\CN^{\circ}$-lenses. Then $$
K^\CN (Z )
\overset{(i_\alpha^*)}\to{\underset{\Sigma i_{\alpha !}}
\to\rightleftarrows}K^\CN (\{ Z_{\alpha}\}):= \mathop\prod\limits_\alpha K^\CN  (Z_\alpha  )
 \tag 3.4.1$$ are homotopically mutually inverse weak equivalences of spectra; 
 the arrow in (3.4.1) is transitive with respect to the refinement of
partitions. 
\endproclaim

{\it Proof.} The argument repeats that from 3.2 with $\CP(\cdot )$ replaced by $\CP (\cdot)^{\CN}
$,  using 2.7(ii) to refine
$\{ Z_{\alpha }\}$ to a partition that comes from a closed filtration. \hfill$\square$

\medskip

{\bf 3.5.} {\it From now on,  $X$ is compact.} Let $U$ be an open subset of $X$, and 
suppose that we have a round cone  $\CN =\CN_{W}\subset TW$ defined over an open neighborhood
$W$ of $X\smallsetminus U$. The complexes  $F$ with $SS(F)\cap$Int$(\CN^{\vee} )=\emptyset$
 form a thick subcategory $\CP (X)^{\CN} \subset\CP (X)$;
 set $K^\CN (X):= K(\CP (X)^{\CN})$. Replacing $X$ by $U$ and $\CN$ by $\CN_{W\cap U}$, we get $
 \CP (U)^{\CN}$
 and $K^\CN (U)$.
 
 \proclaim{\quad Theorem-construction} There is a natural homotopy morphism of spectra $\varepsilon^\CN_{ U}:
 K^\CN (U)\to K$, defined by the datum of  $(U,\CN_{U\cap W})$, such that the diagram  
$$\spreadmatrixlines{2\jot}
\matrix
K^\CN (X)& \lra
 &
K(X) \\  
  j^{*}_{U}\downarrow\,\, \,\,\,\,&&\,\,\,\,\,\,\,\,\,\downarrow   R\Gamma \\
K^\CN (U)  &
    \buildrel{\varepsilon^\CN_{ U}}\over\lra
 &
K 
\endmatrix
\tag 3.5.1$$
is naturally homotopically commutative.
\endproclaim

  {\it Proof.} We construct $\varepsilon^{\CN}_U$ and the homotopy in (3.5.1)
using an auxiliary datum of $P$, $\{ Z_\alpha \}$ below, and then show that the output is independent of the choice.

Choose an $\CN^\circ$-constructible subset  $P\subset W$ which is 
a neighborhood of $X\smallsetminus U$.  Choose a locally closed subanalytic partition $\{ Z_\alpha \}$ of $X$ which is finer than the partition $\{ P, X\smallsetminus P\}$. Let $\{ Z_\alpha^{in}\}$, $\{ Z_\alpha^{out}\}$ be the components that lie in $X\smallsetminus P=U\smallsetminus P$, resp.~in $P$.
Assume that each $Z_\alpha^{out}$ is an $\CN^\circ$-lens. Such a $\{ Z_\alpha \}$ exists by 2.8. We write $K(\{ Z_\alpha^{in}\}):= \prod K(Z_\alpha^{in})$, etc.

Let $\varepsilon^{\CN }_{U}$
be the composition $K^\CN (U)\to K(U)\buildrel{i^{*}_{\alpha}}\over\lra K(\{ Z_\alpha^{in}\}) \buildrel{R\Gamma_{c}}\over\lra K$. Here the last arrow is the map with components $R\Gamma_c : K(Z_\alpha^{in})\to K$.

Let us construct the homotopy in (3.5.1). By Remark in 3.2,  $R\Gamma : K(X)\to K$ is naturally homotopic to the composition
$K(X)\buildrel{(i^*_\alpha )}\over\lra K (\{ Z_\alpha^{in}\})\times K(\{ Z^{out}_\alpha \})\buildrel{R\Gamma_c}\over\lra K$. Therefore the composition $K^\CN (X)  \to K(X)\buildrel{R\Gamma}\over\lra K$ is the sum of $\varepsilon^{\CN}_{ U} j^*_U$ and the map
defined by the DG functor $\CP(X)^\CN \to \CP$, $F\mapsto \Sigma\, R\Gamma_c (Z_\alpha^{out}, i^*_\alpha F)$. Since, by 2.5, the latter functor take values in acyclic complexes, 
the map $K^\CN (X) \to K$ it defines is naturally homotopic to zero, and we are done.

It remains to show that the construction does not depend on the choice of $(P, \{ Z_\alpha \})$. The latter datum form an ordered set (one has $(P, \{ Z_\alpha \})\le (Q, \{ Z_\beta \})$ if $P\supset Q$ and $\{ Z_\alpha \}\le \{ Z_\beta \}$) whose nerve is contractible.
 So it suffices to identify the morphisms $\varepsilon^\CN_{ U}$ and the homotopies in (3.5.1) constructed by means of $(P, \{ Z_\alpha \})\le (Q, \{ Z_\beta \})$ in a transitive manner. Refining $\{ Z_\alpha \}$ and $\{ Z_\beta \}$, we can assume that $\{ Z_\alpha \}$ is finer than the partition $\{ Q, P\smallsetminus Q, X\smallsetminus P \}$. Then the maps $\varepsilon^\CN_{ U}$ for $(P,\{ Z_\alpha \})$ and $(Q, \{ Z_\beta \})$ differ then by a map $K^\CN (X) \to K$ defined by the DG  functor $\CP (X)^\CN \to \CP$ which assigns to $F$ the direct sum of
$R\Gamma_c (Z_\alpha , i^*_\alpha F)$ for all $\alpha$ with $Z_\alpha \subset P\smallsetminus Q$. As above,  this functor takes values in acyclic complexes, so the map of the $K$-spectra is naturally homotopic to zero. The homotopies in (3.5.1) are identified in the same way.
 \hfill$\square$

\medskip

{\bf 3.6.} {\it Remarks.} (i)   By the construction, if 
$U$ is the disjoint union of finitely many $ U_{i}$'s,
then $K^\CN (U)=\oplus K^\CN (U_{i})$ and the $i$th component of $\varepsilon^\CN_{ U}$ equals $
\varepsilon^\CN_{ U_{i}}$.

(ii) Every complex $F\in\CP (X)^{\CN}$ yields a homotopy point $[F]\in K^\CN (X)$, hence a homotopy
point $\varepsilon^{\CN}_{ U} ( F) :=\varepsilon^{\CN}_{ U}   (j^*_{U}[G])\in K$. The homotopy from 
 (3.5.1) yields a natural identification of  the homotopy points $$[R\Gamma (X,F)]=  \varepsilon^{\CN}_{ U} ( F)
  \in K. \tag 3.6.1$$ Therefore if $U=\sqcup U_{i}$ as in (i), then $[R\Gamma (X,F)]= \Sigma\,
   \varepsilon^{\CN}_{ U_{i}} ( F)$.
   
   (iii)  The construction is natural with respect to embeddings of the cones $\CN$, and shrinking of $U$ and $W$
   (subject to the condition $U\cup W=X$). 
   
 (iv) If $R$ is commutative and $\CP$ consists of perfect complexes, then we have the morphism of the Picard groupoids $\det : \Pi K(R)\to \CL (R)$ (see 1.8(vii)), hence a graded super line $\CE^{\CN }_{U}(F):=\det \varepsilon^{\CN}_{ U} ( F)$ identified with $\det R\Gamma (X,F)$. If $U=\sqcup U_i$, then, by (i), one has a canonical factorization isomorphism $\CE^{\CN}_{ U}(F) =\otimes \CE^{\CN}_{ U_i}(F)$, hence
 we get the $\varepsilon$-factorization $\det R\Gamma (X,F)\iso \otimes \CE^{\CN}_{ U_i}(F)$.

(v) The construction of 3.5 should be compatible with the Poincar\'e-Verdier duality. Namely, all
 the constructions  can be literally repeated replacing the datum of functors $(i_! ,i^* , R\Gamma_c )$ by the dual datum $(Ri_* , Ri^! , R\Gamma)$, and replacing
 $\CN^{\circ}$-lences  by $\CN$-lenses. Presumably, the resulting $\varepsilon$ morphism and the
compatibility homotopy of (3.5.1) do not change; the key point should be the lemma in 2.5. If $\CP$ consists of perfect complexes, then  the
Poincar\'e-Verdier duality $\CP_{R}(X)^{\circ}\iso \CP_{R^{\circ}}(X)$,\footnote{Here $R^{\circ}$ is
$R$ with the opposite product.}  interchanges $(i_! ,i^* , R\Gamma_c )$ and $(Ri_* , Ri^! , R\Gamma)$ and, by 2.5, identifies $(\CP_{R }(X)^{\CN })^{\circ}$ 
 with $\CP_{R^{\circ}}(X)^{\CN^{\circ}}$. The should imply that
 the  duality identifies  $\varepsilon^{\CN }_{U}$ 
 with $\varepsilon^{\CN^{\circ }}_{U}$, and interchanges the respective homotopies in (3.5.1).  I did not check the details.

\medskip

{\bf 3.7.} {\it Variant.}   Let $Y\subset X$ be a closed subset, and $\nu$ be a continuous nowhere vanishing 1-form on $W:=X\smallsetminus Y$. We have the category $\CP (X)^{\nu }$ of $ \CP$-constructible complexes $F$ on $X$ such that $SS (F)\cap\nu (W)=\emptyset$ 
 and its  $K$-spectrum $K^\nu (X)$. There is a similar category $\CP (U)^{\nu }$ for each open neighborhood $U$ of $Y$. Let $\CP (\tilde{Y})^{\nu}$ be the inductive limit of the directed set of these categories and   $K^\nu (\tilde{Y})$ be  its $K$-spectrum (which is the colimit of $K^\nu (U)$'s). One has  evident morphisms $K^\nu (X)\to K(X)$ and $j_Y^* : K^\nu (X) \to K^\nu (\tilde{Y})$. 

\proclaim{\quad Theorem} There is a natural homotopy morphism $\varepsilon_{\nu Y}: K^\nu (\tilde{Y}) \to K$
such that the diagram  
$$\spreadmatrixlines{2\jot}
\matrix
K^\nu (X)& \lra
 &
K(X) \\  
  j^{*}_{Y}\downarrow\,\, \,\,\,\,&&\,\,\,\,\,\,\,\,\,\downarrow   R\Gamma \\
K^\nu (\tilde{Y})  &
    \buildrel{\varepsilon_{\nu Y}}\over\lra
 &
K 
\endmatrix
\tag 3.7.1$$
is naturally homotopically commutative.
\endproclaim

  {\it Sketch of a proof.} Consider the set $\CC =\CC (\nu )$ of all round cones $\CN \subset TW$ such that $\nu\in$ Int$(\CN^\vee )$; the map $\CN \mapsto \CN^\vee$ identifies $\CC$ with the set of all round cone neighborhoods of $\nu (W)$. Our $\CC $ is ordered by inclusion; it is directed, and $\CP (X)^\nu =\cup_\CC \CP (X)^\CN$. To establish the theorem, one combines the constructions from 3.5 for all $\CN\in\CC$. Namely, for  an open neighborhood $U$ of $Y$ we consider the set $\CC_U$  of pairs $(\CN ,P)$ where $\CN\in\CC$ and $P\subset W$ is an $\CN$-constructible subset which is a neighborhood of $X\smallsetminus U$. Then $\CC_U$ is ordered by inclusion, and is  directed. For each $(\CN ,P)$ we consider the morphism $\varepsilon^{\CN }_{U}$ and the homotopy from (3.5.1). If $(\CN ,P)\le (\CN' ,P')$ then the restriction of $\varepsilon^{\CN' }_{P' }$ to $K^\CN (U)$ is naturally homotopic to $\varepsilon^{\CN}_{ U}$, same for the homotopies from (3.5.1). The identifications are transitive, which yields the pull-back of $\varepsilon_{\nu Y}$ to $K^\nu (U)$ and the corresponding homotopy. Their compatibility for different $U$'s is immediate. The details are left to the reader.   \hfill$\square$
    
   \medskip

{\bf 3.8.} We are in the setting of  3.7. Suppose that there is a $C^1$-function $f$ on $X$ which is locally constant on $Y$ and such that $\nu =df|_W$ (the principal case is that of finite $Y$).
 Let us show that in this situation the datum of 3.7 can be recovered from 
 the classical Morse theory picture. The latter can be summarized as follows:
 
 Let $a_{0}>\ldots >a_{n}$ be the  values of $f$ on $ Y$; set $Y_i := Y_{f=a_i}$. Any
 $F\in D(X)$ carries a canonical filtration $F_{0}
  \subset\ldots\subset F_{2n}=F$ with $F_{2i-1}\subset F_{2i}\subset F$ being the Morse
  filtration for $U=X_{f>a_{i}}$, $V=X_{f<a_{i}}$ (see 2.1). If
   $F\in \CP (X)^{\nu }$, then, by 2.5, $\gr_{2i}F=M_{f=a_{i}}(F)$ is supported on $Y_i$. Set  $M_{i}(F):= R\Gamma (X,  M_{f=a_{i}}(F))$.
By 2.5 and the lemma from 2.1,
 one has $R\Gamma (X,\gr_{2i+1}F)=0$.  Therefore we get the {\it Morse filtration} on
 $R\Gamma (X,F)$ with 
 $\gr_i R\Gamma (X,F)= M_{i}(F)$. By (i) of the lemma below, one has
 DG functors $M_{i}: \CP (X)^{\nu }\to \CP$, hence
   the   morphisms of the $K$-spectra $M_{i}:K(X)^{\nu }\to K$; the Morse filtration on
  $R\Gamma$ yields a {\it Morse identification} of the homotopy morphisms
   $R\Gamma =\Sigma\, M_{i}: K(X)^{\nu }\to K$. Of course, $M_{i}$
   is defined on  $\CP (U_{i })^{\nu }$ where $U_{i}$ is any open neighborhood of $Y_i$,
   hence we have $M_{i}:K(U_{i})^{\nu }\to K$.
   
     \proclaim{\quad Proposition} (i) 
   For $F$ in $ \CP (X)^{\nu }$ or $\CP (U_i )^{\nu }$ one has
    $M_{i}(F)\in\CP$.
   
   (ii) The above morphisms of the $K$-spectra  
    $M_{i}$ identify naturally with $\varepsilon_{\nu Y_i} $, and the 
   Morse identification $R\Gamma =\Sigma\, M_{i}$ with the homotopy from (3.7.1).
    \endproclaim
    
    {\it Proof.}    Let $U_i$ be 
     any open neighborhoods  of $Y_i$ such  that $U_i$ does not
  intersect $U_j$ and  $X_{f=a_j}$ for $j\neq i$; set $U:=\sqcup U_{i}$. Pick $\CN\in\CC (\nu )$ (see the proof 3.7); we will construct the identifications of (ii) on $\CP (X)^\CN$ leaving the compatibilities for different $\CN$'s to the reader.

   We will construct a  filtration $  X_0 \subset X_1\subset \ldots \subset X_{2n}=X$ such that $X_i$ are open subanalytic subsets, the closure of $Z_i := X_{2i} \smallsetminus X_{2i-1}$ lies in $U_i$, Int$(Z_i )\supset Y_i$, and $X_i \cap W$ are $\CN^\circ$-open subsets of $W$. Such a datum yields then (i) and (ii).\footnote{We leave to the reader to check that the identifications of (ii) constructed from different filtrations $X_i$ are naturally homotopic.} Indeed, 
    on any sheaf $F$  we get a filtration $F'_0 \subset \ldots \subset F'_{2n}=F $, $F'_i := j_{X_i !} j_{X_i}^* F$. If $F\in\CP (X)^\CN$, then, by the lemma in 2.5, one has $\gr'_\cdot F \in \CP (X)^\CN$. It is clear that $M_i (\gr'_a F )$ equals $M_i (F)$ if
   $a=2i$ and vanishes otherwise. Thus $M_i (F)=  R\Gamma (X, \gr'_{2i}F)\in \CP$, hence (i), and $R\Gamma (X,\gr'_{2i+1}F)=0$. The filtration of $R\Gamma (X,F)$ by $R\Gamma (X, F'_{2i})$ coincides with the Morse filtration in the filtered derived category.\footnote{
 Since $\gr_i \gr'_j R\Gamma (X,F) $ equals $M_i (F)$ if $i=j$ and vanishes otherwise.} By the construction in 3.5, $\gr'_i R\Gamma (X,F)= R\Gamma_c (Z_i , i^*_{Z_i}F)=\varepsilon^{\CN}_{ U_i }(F)$, and the homotopy of (3.5.1) comes from our filtration on $R\Gamma (X,F)$, and we are done. 
   
 To construct $X_i$, pick any 
   non-negative $C^1$-function $\phi$ on $X$ which equals 1 on a neighborhood of $Y$ and 0 on a neighborhood of $X\smallsetminus U$. For a constant $\epsilon >0$ set $f^\pm_\epsilon := f\pm \epsilon \phi$. If $\epsilon$ is sufficiently small, then
   $df^\pm_\epsilon (W)\subset $ Int$(\CN )$. Consider a filtration $X'_0 \subset \ldots X'_{2n} =X$ where $X'_{2i}:= X_{f_\epsilon^+ > a_{i}}$, $X'_{2i -1}:=  X_{f_\epsilon^- > a_i} $. This filtration satisfies all our conditions except subanaliticity. We modify it  as follows. For  $x\in X_{f^+_\epsilon = a_i}$  pick, as in the proof of 2.7(i), a subanalytic $\CN^\circ$-lens $P_x^{+}\subset W$ such that $x\in $ Int$(P_x^+ )$ and $ f^+_\epsilon (\bar{P}_x^+ \smallsetminus P^+_x )> a_i$; we demand that $P^+_x$ does not intersect $U_j$ and $X_{f=a_j}$ for $j\neq i$.     Replacing $+$ by $-$, pick
similar $P_x^-$ for $x$ in the closure  $\Gamma_i$ of $ X_{f^-_\epsilon = a_i} \smallsetminus X_{f^+_\epsilon = a_i}$; we want in addition that $P_x^- \subset U_i$. Let $\{ P^{+}_{i\alpha }\}$ be a finite subset of $P^+_x$'s that cover $X_{f^+_\epsilon = a_i}$, and $\{ P^-_{i\beta} \}$ that of $P^-_x$'s that cover $\Gamma_i$. Then $X_{2i}:= X'_{2i}\smallsetminus \cup P^+_{i\alpha}$, $X_{2i-1}:= X_{2i}\smallsetminus \cup P^-_{i\beta}$ is the promised filtration.
   \hfill$\square$

 \bigskip

\centerline {\bf 4. The animation of the Dubson-Kashiwara formula.}

\medskip

In this section we explain 0.1(b). To do this, we first ``localize" the map of $K$-spectra $R\Gamma : K(X)\to K$ defining the map $\varepsilon$, which can be further microlocalized.

{\bf 4.1.} Our $X$ is a compact real analytic manifold. We use the notation of 1.3.

\proclaim{\quad Proposition-construction} There is a natural homotopy morphism 
of spectra $$\varepsilon  : K (X )\to C_\n (X, K ) \tag 4.1.1 $$ such that the composition
$K(X)\buildrel{\varepsilon}\over\to C_{\n}(X,K)\buildrel{tr}\over\to K$ is naturally
homotopic to
$R\Gamma$.
\endproclaim

 {\it Proof.} 
We construct $\varepsilon$ using an auxiliary datum of $\{ Z_\alpha , V_\alpha \}$
 where $\{ Z_\alpha \} $ is a locally closed subanalytic partition of $X$ and each 
 $V_\alpha$ is a contractible subset of $X$ which contains $Z_\alpha$.  
 The set $\CA (X)$ of such data is ordered with respect to refinement:  $\{ 
Z_\alpha , V_\alpha \} \le \{ Z_\beta , V_\beta \}$ if for every $\beta$ there is a (necessary unique) $\alpha=\alpha (\beta )$ such that $Z_\beta \subset Z_\alpha$ and $V_\beta \subset V_\alpha$. 

Our $\CA (X)$ is directed. Indeed, let Par$^0 (X)\subset $ Par$^a (X)$ be 
the subset  of locally closed partitions with each $Z_\alpha$ contractible. 
One has an embedding  of ordered sets Par$^0 (X) \hra  \CA (X)$, 
$\{ Z_\alpha \} \mapsto \{ Z_\alpha , Z_\alpha \}$. Then Par$^0 (X)$ is directed (see 3.1)
and is a cofinal subset of $ \CA (X)$, whence the assertion. 

Consider a chain of morphisms of spectra
$$ K (X)\buildrel{\Sigma i_{\alpha !}}
\over\longleftarrow  
K (\{ Z_\alpha \} )\buildrel{(R\Gamma_c )}\over\lra \mathop\prod\limits_\alpha K
 \buildrel{(tr)}\over\longleftarrow \mathop\prod\limits_\alpha C_\n 
 (V_\alpha , K )\to C_\n (X,K ), \tag 4.1.2$$ where the middle arrows are diagonal matrices,
last arrow comes from the embeddings $V_\alpha \hra X$.\footnote{As in 3.1, 
$R\Gamma_c$ is constructed using a natural flabby resolution $F\to \tilde{F}$.} 
Both arrows that look to the left are weak equivalences (the first one by 
 3.2, the second one since $V_\alpha$ are contractible), so they admit canonical homotopy inverses.
Let $\varepsilon_{\{ Z_\alpha ,V_\alpha \}} : K (X) \to C_\n (X, K )$ be the 
 composition. 

For $\{ Z_\beta , V_\beta \}\ge \{ Z_\alpha , V_\alpha \} $ there is a natural 
 morphism of  diagrams  $$\kappa : (4.1.2)_{\{Z_\beta , V_\beta \}}\to 
 (4.1.2)_{\{Z_\alpha ,V_\alpha  \}} \tag 4.1.3$$ which is identity at the ends $K (X)$ and
 $C_\n (X,K )$, and is formed by morphisms with the components 
 $i_{\beta \alpha (\beta ) !} : K (Z_\beta )\to K (Z_{\alpha (\beta )})$, $ id_{K} : 
 K^{ \beta} \to K^{ \alpha (\beta )}$,  $C_\n (V_\beta , K )\to 
 C_\n (V_{\alpha (\beta )}, K )$ in the middle.  This $\kappa$  yields a natural homotopy 
 between 
$\varepsilon_{\{ Z_\alpha ,V_\alpha \}}$ and $\varepsilon_{\{ Z_\beta ,V_\beta \}}$.
 Our $\kappa$'s are  transitive, so  (4.1.2)$_{\{ Z_\alpha , V_\alpha \}}$ form an
  $\CA (X)$-diagram. Since $\CA (X)$ is directed, we see that 
  all $\varepsilon_{\{ Z_\alpha ,V_\alpha \}}$'s for $\{ Z_\alpha ,V_\alpha \}\in\CA (X)$
   are all naturally homotopic, thus defining (4.1.1). 

 A slightly more canonical way to define $\varepsilon$ (which will be useful in 4.6)
  is to consider a chain of morphisms of spectra
  holim$_{\CA (X)}$ (4.1.2)$_{\{ Z_\alpha , 
 V_\alpha \}}$. Add  the standard  weak equivalences 
 $K (X) \to $ holim$_{\CA (X)} K (X)$ at the left and 
 holim$_{\CA (X)}C_\n (X,K )\leftarrow C_\n (X,K )$ at the right, and invert homotopically the arrows 
 looking to the left. Our $\varepsilon$ is the composition.

It remains to identify the composition $tr\,\varepsilon$ with $R\Gamma$. 
Notice that each term of (4.1.2) projects
naturally to $K$: the projection is $R\Gamma$ at the left end, $tr$ at the right end, and its $\alpha
$-components at the middle terms are, respectively, $R\Gamma_{c}$, $id_{K}$, and $tr$. The projections
commute with the arrows in (4.1.2), so they yield the identification $tr\,\varepsilon_{\{ Z_{\alpha},V_{\alpha}\}}=R\Gamma$.
Since they are also compatible with the morphisms from (4.1.3), the choice of 
$\{ Z_\alpha ,V_\alpha \}$ is irrelevant.  \hfill$\square$

 {\it Remarks.} (i) One can define $\varepsilon$ using {\it any} directed subset of $\CA (X)$,
 e.g.~Par$^0 (X)$ (which has advantage of excluding the redundant 
 $V_\alpha$'s),  or the
 subset of all $\{ (Z_{\alpha },V_{\alpha})\}$ with open $V_{\alpha}$'s.

(ii) Realizing the homotopy inverse to $\Sigma i_{\alpha !}$ as  $(i^*_\alpha )$ (see 3.2), 
one can write $\varepsilon_{\{ Z_\alpha ,V_\alpha \}}
 =\mathop\sum\limits_\alpha \varepsilon_{( Z_\alpha ,V_\alpha )}$ where each
  $\varepsilon_{( Z_\alpha ,V_\alpha )}$ is the individual 
  composition $$ K (X)\buildrel{ i^*_{\alpha }}\over\lra  K 
  (Z_\alpha )\buildrel{R\Gamma_c}\over\lra  K \buildrel{tr}
  \over\longleftarrow  C_\n (V_\alpha , K )\to C_\n (X,K ). \tag 4.1.4$$
 We did not do this on the spot for $\kappa$ of (4.1.3) {\it strictly} commutes 
 with all the arrows of (4.1.2) as written. A choice of a point in $V_\alpha$ provides a morphism   $K\to C_\n (V_\alpha , K )$ which is right inverse, hence homotopy inverse, to the left arrow $tr$ in (4.1.4).
 
 (iii) Any $F\in\CP (X)$ yields a homotopy point $[F]\in K(X)$, hence a homotopy point
 $\varepsilon (F)\in C_{\n}(X,K)$. Truncating $K(X)$ to $K_{0  }:=\pi_{0}K(X)$, it becomes
 a mere zero cycle of degree equal to the Euler characteristics, which is canonically defined up
 to all higher homotopies.
 
 {\it Exercises} (irrelevant for the sequel). (i) Replacing $X$ by a locally closed subset $Z$,  define morphisms $\varepsilon_{!},\varepsilon_{*}
 :K(Z)\to C_{\n}(Z,K)$ with $tr\,\varepsilon_{!}=R\Gamma_{c}$, $tr\,\varepsilon_{*}=R\Gamma$.
Check that $\varepsilon_{!}$ is compatible with $Rf_{!}$ maps, and 
$\varepsilon_{*}$ - with $Rf_{*}$ maps (see 3.1).

(ii) The trace map  $tr  : C_\n (X, K (X))\to K (X)$ admits a 
{\it canonical} homotopy section\footnote{Consider the trace map   $tr :C_\n (X,P)\to P$
 for any spectrum $P$.  Then any point $x\in X$ yields a section $s_x : P\to C_\n (X,P)$ of
  $tr$ which depends on $x$ unless $X$ is contractible.} $s : K (X)\to C_\n (X, K (X)) $ 
  defined as the composition $K (X) 
  \buildrel{(i_{\alpha !} i_\alpha^* )}\over\lra \mathop\prod\limits_\alpha 
  K (X) \buildrel{\sim}\over\leftarrow \mathop\prod\limits_\alpha C_\n (V_\alpha
  , K (X))\to C_\n (X, K (X))$, and $\varepsilon$ is naturally homotopic to the composition 
  $K (X)\buildrel{s}\over\lra C_\n (X, K (X)) \buildrel{R\Gamma}\over\lra C_\n (X,K )$.

\medskip

{\bf 4.2.}  Below we consider the conical topology on $T^* X$, so an open subset $E\subset T^* X$ is always assumed to be conic. Set $U_E := \pi_{T^*}(E)$. For a spectrum $R$ we denote by $R_{T^* X}$ the corresponding constant presheaf on $T^* X$.

One has a cofibrant presheaf of spectra $K^!_X  $ on $X$, 
$K^!_X (U) =C_\n (X, X\smallsetminus U; K )$ (see 1.7). Consider its 
pull-back $\pi^*_{T^*}K^!_X $ to $T^* X$, $\pi^*_{T^*}K^!_X (E)= K^!_X (U_E ) $; one has an evident morphism of presheaves  $C_\n (X,K)_{T^* X} \to   \pi^*_{T^*}K^!_X $. The  morphisms 
$C_\n (X,K ) \to R\Gamma (X, K^!_X )\to  R\Gamma (T^* X,\pi^*_{T^*}
 K^!_X )$ are weak equivalences of spectra (the first one by 1.7, 
 the second one since the fibers of $\pi_{T^*}$ are contractible). 

 For  $E$ as above, let $\CP (X)_{E}^{\dagger}$ be a thick subcategory of $\CP (X)$
 formed by those complexes $F$ that $SS(F)\cap E=\emptyset$. Since
   $\CP (X)_{E}^{\dagger}\subset\CP (X)_{E'}^{\dagger}$ if $E\supset E'$, our  $\CP (X)_{E}^{\dagger}$
   form a presheaf  $\CP (X)^{\dagger}$ of DG categories over $T^{*}X$. Denote by 
   $K^\dagger $ the presheaf of  $K$-spectra $
    K^{\dagger}(E) :=K(\CP(X)^{\dagger}_{E})$; one has an evident morphism   $\iota : K^{\dagger}    \to K(X)_{T^* X}$. We define a morphism  $\eta$ by the next commutative diagram:
   $$ \spreadmatrixlines{2\jot}
\matrix  K(X)_{T^* X}& \buildrel{\varepsilon}\over\lra  & C_\n (X,K)_{T^* X} \\
\uparrow \iota && \downarrow \\
  K^\dagger &  \buildrel{\eta}\over\lra & \pi^*_{T^* } K^!_X .
  \endmatrix \tag 4.2.1$$
Let $\pi^*_{T^* } K^!_X \to F^\mu$ be a fibrant resolution of $\pi^*_{T^* } K^!_X$ for the Jardine ``sheafified" model structure (see 1.6).

\proclaim{\quad Theorem-construction} The composition $K^{\dagger}\buildrel{\eta}\over\to  \pi^*_{T^* } K^!_X \to F^\mu$ is naturally homotopic to zero.  More precisely, it  factors naturally through a homotopically trivial presheaf.  \endproclaim

Here ``homotopically trivial" means triviality of the presheaves of the homotopy groups. 
The proof occupies 4.3--4.6. We assume it for the moment.

 Consider a presheaf  $K^\mu$  on $T^* X$, $ K^\mu (E):=  \CC one (K^\dagger (E)\to K (X))$. It  can be interpreted as follows.
Let $\CP^\mu (E)$ be the Verdier quotient $\CP (X)/\CP (X)^\dagger_E$ (for a construction on DG level, see \cite{Dr}).  This is the category of microlocal constructible sheaves on $E$ as defined in \cite{KS} 6.1.  According to 1.8(vi), one has a canonical weak equivalence of spectra $K^\mu (E)\iso K(\CP^\mu (E))$. 
   
Consider the composition of  $C_\n (X,K)_{T^* X}\to \pi^*_{T^* } K^!_X \to F^\mu$. Let us factor it as $C_\n (X,K)_{T^* X}\to C_\n (X,K)^\flat_{T^* X}\to F^\mu$ where the first arrow is an objectwise weak equivalence, the second one is a fibration for the plain  (non-sheafified) model structure on presheaves. Let $F^\dagger $ be the fiber  of this  fibration. By the above theorem, the composition $K^\dagger  \to K(X)_{T^* X} \buildrel{\varepsilon}\over\lra C_\natural (X,K)^\flat_{T^* X} \to F^\mu$ is naturally homotopic to zero. This 
 homotopy  yields a naturally homotopy commutative diagram  
 $$\spreadmatrixlines{2\jot}
\matrix
K^\dagger  & \buildrel{\varepsilon^\dagger}\over\lra
 &F^\dagger \\  
\downarrow &&\downarrow \\
K (X)_{T^* X} &
    \buildrel{\varepsilon }\over\lra
 &C_\n (X,K)^\flat_{T^* X}\\
\downarrow &&\downarrow \\
K^\mu &
    \buildrel{\varepsilon^\mu }\over\lra
 &F^\mu .
\endmatrix
\tag 4.2.2$$
Notice that $F^\mu (E)=R\Gamma (E, \pi^*_{T^* }K^!_X  )$ and $F^\dagger (E)=
 R\Gamma_S (T^* X, \pi^*_{T^* }K^!_X  )$ where 
$S:= T^* X\smallsetminus E$. The morphism $\varepsilon^\mu$ 
 can be seen as   
 the microlocalization of the $\varepsilon$ map.

{\it Remarks.} (i) I do not know if $\CP^\mu (E)$ is always Karoubian, or if the 
Karoubian property holds locally on $T^* X$. If this is not the case, it would be nice  to extend $\varepsilon^{\mu}$ to the
 presheaf of $K$-spectra of the idempotent completions.

(ii) The theorem-construction from 3.5 (and 3.7) is a corollary of the present theorem.\footnote{The  constructions from 3.5 are presented independently for they are less involved and more direct than those used in the proof of 4.2.} Indeed,
in the notation of loc.~cit., set $E:=$ Int$(\CN^{\vee})\subset T^{*}X$. Then $\CP (X)^\CN =\CP (X)^\dagger_E$,
$K^{\CN}(X)=K^\dagger (E)$ and one has $C(X,X\smallsetminus W;K)=R\Gamma (W,K^{!}_{X})\iso R\Gamma (E,\pi^{*}_{T^{*}}K^{!}_{X})$
 (since the $E/W$ is fiberwise contractible).
 Therefore the composition $K^{\CN}(X) \to K(X)
 \buildrel{\varepsilon}\over\lra
 C_{\n}(X,K)\to C_{\n}(X, X\smallsetminus W;K) $ is naturally homotopic
 to zero by the theorem. Since $ C_{\n}(X\smallsetminus W,K)$ is the homotopy fiber of the last arrow,
 we see,
 as in (4.2.2), that the homotopy amounts to a lifting of $ \varepsilon|_{K^\CN (X)}$
 to a morphism $ K^\CN (X)\to C_{\n}(X\smallsetminus W,K)$. We leave it to the reader to 
 check that the composition $ K^\CN (X)\to C_{\n}(X\smallsetminus W,K)\buildrel{tr}\over\lra
 K$ is naturally homotopic to $\varepsilon^{\CN}_{ U}j^{*}_{U}$ of (3.5.1), to
 construct $\varepsilon^{\CN}_{ U}$ itself, and to show that the homotopy from (3.5.1) comes from
 the identification $R\Gamma = tr\,\varepsilon$ of 4.1.

(iii) Let $F$ be any $\CP$-constructible complex. It yields a homotopy point $[F]\in K^\dagger (T^* X\smallsetminus  SS(F))$, hence
a homotopy point   $\varepsilon^\dagger (F)\in R\Gamma_{SS(F)}(T^* X, \pi^*_{T^*} K^!_X )$ whose image by  $R\Gamma_{SS(F)}(T^* X, \pi^*_{T^*} K^!_X )\to R\Gamma (T^* X, \pi^*_{T^*} K^!_X )\iso R\Gamma (X, K^!_X )\buildrel{\sim}\over\leftarrow C_\n (X,K)$ equals $\varepsilon (F)$; here the middle arrow is pull-back by the zero section $X\to T^* X$. One can view  $\varepsilon^\dagger (F)$ as an animation of Kashiwara's characteristic cycle $CC(F)$; therefore the identification $[R\Gamma (X,F)]=tr\, \varepsilon (F)$ of 4.1 becomes an animation of the Dubson-Kashiwara formula (0.1.3).  To see this, notice that for any abelian group $A$ the group $H^0 R\Gamma_{SS(F)}(T^* X, \pi^*_{T^*} A^!_X )$ is 
the group of Lagrangian cycles supported on $SS(F)$ with coefficients in $A$ (see \cite{KS} 9.3). Now the image of $\varepsilon^\dagger (F)$ by the map $K\to \pi_0 (K)=K_0 (P)$ coincides with $CC(F)$.
This follows easily from (ii) and the Morse-theoretic interpretation of $\varepsilon_{df}$ from 3.8 (for $\nu =df$ intersecting transversally a component of $SS(F)$). The details are left to the reader.

(iv) When $R$ is commutative and $\CP$ consists of perfect complexes,  the image of $\varepsilon^\dagger (F)$ by the determinant map provides a microlocal description of the determinant line $\det R\Gamma (X,F)$. Can one deduce from it  the Lefschetz formulas  of \cite{KS} 9.6?

{\bf 4.3.} {\it Proof of the theorem.}    We will use the next 
approximation of the topology of $X$ and the conical topology of $T^* X$:

Let $X^a$ be the Grothendieck topology formed by open subanalytic subsets of $X$ 
(the morphisms are  embeddings, the coverings are evident ones). 

Let $T^* X^a$ be the Grothendieck topology whose objects are
 open   $D\subset T^* X$ which are either $\pi_{T^*}$-preimages of subsets 
 from $X^a$ or {\it subanalytic lenses} as defined below. The morphisms are  embeddings $D_{1}\subset
 D_{2}$ subject to a condition: if $D_{1}$ is a lens, $D_{2}$ is not, then
 one has $\bar{U}_{D_{1}}\subset U_{D_{2}}$. The coverings are evident ones.

 We say that $D$ is a {\it lens}
 if there is a round cone $\CN^{(D)}$ defined over an open neighborhood $W$ of $\bar{U}_{D}$
 and an $\CN^{(D)}$-lens $Z$  such that $U_{D}=$ Int$(Z)$, 
  $D=$ Int$(\CN^{(D)\vee}_{U_{D}})$.  By 2.4(i), $Z$ is uniquely determined by $D$; we write
 $U_{D}^{+} :=Z$, $U_{D}^{-}:=Z^{\circ}$, and refer to $U_{D}$ and $U^{\pm}_{D}$ as
  the {\it open} and {\it locally closed
 $D$-lenses.} We say that $D$ is a subanalytic lens if $Z$ is a subanalytic $\CN^{{(D)}}$-lens (see 3.4).
Notice that $\bar{D}= \CN^{(D)}_{\bar{U}_{D}}$, and one can take for $\CN^{(D)}$ any extension
of $D^{\vee}$ to a neighborhood of $\bar{U}_{D}$ (see 2.4(ii)).

For any open $U\subset X$ we denote by $U^a$ the topology of open subsets $V$ of $U$ that belong to $X^a$. For a presheaf of spectra $M$ on $X$ we define $M^a (U) $ as holim$\, M(V)$, $V\in U^a$. Then $M^a$ is a presheaf of spectra on $X$. One has an evident natural morphism $\alpha: M\to M^a$. If $M$ is fibrant for the ``sheafified" model structure, then such is $M^a$, and $\alpha$ is an objectwise weak equivalence (to see this, notice that
 every point  admits a base of  neighborhoods from $X^a$, so every hypercovering 
of $U$ can be refined to an $U^a$-hypercovering).

Similarly, for an open  $E\subset T^* X^a$ its open subsets $D$ from $T^* X^a$ form a topology $E^a$. 
For a presheaf $L$ on (the conical topology of) $T^* X$ we get a presheaf $L^a$ where $L^a (E):=  $ holim$\, L(D)$, $D\in E^a$, and
a morphism  $\alpha : L\to L^a$ which is an objectwise weak equivalence if $L$ is fibrant for the ``sheafified" model structure (use 2.4(iii)). 

{\bf 4.4.} If $  \pi^*_{T^* } K^!_X \to F^\mu$ is a fibrant resolution of $\pi^*_{T^* } K^!_X$, then so is $   \pi^*_{T^* } K^!_X \to F^\mu \to F^{\mu a}$. The composition
 $K^{\dagger}\buildrel{\eta}\over\lra \pi^*_{T^* } K^!_X \! \to F^{\mu a}$ can be rewritten as $K^{\dagger}\buildrel{\alpha}\over\lra K^{\dagger a} 
 \buildrel{\eta^{ a}}\over\lra (\pi^*_{T^* } K^!_X )^a \to F^{\mu a}$. Therefore the theorem follows from the next assertion:

\proclaim{\quad Proposition} The morphism
 $\eta^{ a}:  K^{\dagger a} \to ( \pi^*_{T^* } K^!_X )^a $ naturally factors through a homotopically trivial presheaf.
\endproclaim

Notice that the  topology  on $T^* X$  plays  no role in the claim: we deal with mere presheaves.

The proof of the proposition takes the rest of the section.

{\bf 4.5.} We need a lemma.
For $D\in T^{*}X^a$ let $\CA_D \subset \CA (X)$ be the subset of $\{ Z_{\alpha}, V_{\alpha}\}
$ such that $V_{\alpha}$  are open subsets,  
for every $\alpha$ with $Z_\alpha \subset X\smallsetminus U_D$ the set $V_\alpha \smallsetminus
    U_D$ is contractible, and the next condition is satisfied: 
if $D=\pi_{T^{*}}^{-1}(U_D)$, then  the partition $\{ Z_\alpha \}$ is finer than that 
$\{ U_D , X\smallsetminus U_D \}$;  if $D$ is a lens, then   $\{ Z_\alpha \}$
 is finer than $\{ U_{D}^- , X\smallsetminus U_{D}^- \}$, and each $Z_\alpha \subset U_{D}^-$ is an $\CN^{(D)\circ}$-lens.
   For an $n$-simplex $D_\cdot \! := (D_0 \subset \ldots \subset D_n )$ 
   of $ \CN er \, T^{*}X^a$ set $\CA_{D_\cdot}\! := \cap \CA_{D_i}$. 

\proclaim{\quad Lemma}   $\CA_{D_\cdot }$ is a non-empty directed subset of $\CA (X)$.
\endproclaim

{\it Proof.}  There is $a\in [0,n]$ such that $D_{i}$ are lenses for $i\le a$, and are not
 lenses for $i>a$. Set $U_{i}:=U_{D_{i}}$; for $i\le a$
 set $\CN_{i} := \CN^{(D_{i})}$. 
 
To show that   $\CA_{D_{\cdot}}\neq \emptyset$, we present some $\{ Z_{\alpha},  V_{\alpha}\}\in \CA_{D_{\cdot}}$.  Let
  $P_{i}$, $i=0,\ldots , n+1$, be the partition of $X$ defined by the filtration $U^{-}_{0}\subset
  \ldots \subset U^{-}_{a}\subset U_{a+1}\subset \ldots \subset U_{n+1}
  := X$ (see 2.6).
   Choose collections of open subsets $\{ V_{\delta_{0}}\},\ldots , \{ V_{\delta_{n+1}}\}$
   such that the sets $V_{\delta_{i}}$
    and $V_{\delta_{i}}\smallsetminus U_{i-1}$ are contractible, 
and for each $i$ the union $\cup 
V_{\delta_{i}}$ contains $P_{i}$. By 
2.7(ii), for each $i\le a$ one can find a partition $\{ Z_{\alpha_{i}}\}$ of 
$U^{-}_{i}\smallsetminus U^{-}_{i-1}$
by $\CN_{i}^{\circ}$-lenses such that every $Z_{\alpha_{i}}$
 lies in some $V_{\alpha_{i}}$ from 
$\{ V_{\delta_{i}}\}$. For  $i>a$ let $\{ Z_{\alpha_{i}}\}$ be any subanalytic partition of
$U_{i}\smallsetminus U_{i-1}$ if $i>a+1$, or of $U_{i}\smallsetminus 
U^{-}_{i-1}$ if $i=a+1$,
such that each $Z_{\alpha_{i}}$ lies in some 
$V_{\alpha_{i}}$ from $\{ V_{\delta_{i}}\}$. The promised 
$\{ Z_{\alpha},  V_{\alpha}\}$ is the union 
of  all $\{ Z_{\alpha_{i}},  V_{\alpha_{i}}\}$.

Let us check that $\CA_{D_{\cdot}}$ is directed. Pick any $\{ Z_{\alpha},  V_{\alpha}\},\{ 
  Z_{\beta},  V_{\beta}\}\in
  \CA_{D_{\cdot}}$; let us construct their upper bound. Let $Z_{\alpha_{i}}$, etc., be those $Z_{\alpha}$ that lie in $P_{i}$.
   The least upper bound of the partitions is formed by sets $Z_{\alpha_{i}}\cap Z_{\beta_{i}}$,
   which are $\CN_{i}$-lenses for $i\le a$.  Choose open contractible 
    $V(\alpha_{i}\beta_{i})_{\delta}\subset V_{\alpha_{i}}\cap V_{\beta_{i}}$ which cover $Z_{\alpha_{i}}\cap Z_{\beta_{i}}$
  and such that $V(\alpha_{i}\beta_{i})_{\delta}\smallsetminus U_{i-1}$ are also 
  contractible. Choose a subanalytic partition $\{ Z(\alpha_{i}\beta_{i})_{\gamma}\}$ of 
   $Z_{\alpha_{i}}\cap Z_{\beta_{i}}$ such that each $ Z(\alpha_{i}\beta_{i})_{\gamma}$
   lies in some $V(\alpha_{i}\beta_{i})_{\gamma}$; if $i\le a$, then we ask $ 
   Z(\alpha_{i}\beta_{i})_{\gamma}$
   to be $\CN_{i}$-lenses, which is possible by 2.7(i). Now 
   $\{  Z(\alpha_{i}\beta_{i})_{\gamma}, V(\alpha_{i}\beta_{i})_{\gamma}\}
   $ is an upper bound of  $\{ Z_{\alpha},  V_{\alpha}\},\{ 
  Z_{\beta},  V_{\beta}\}$ in
 $ \CA_{D_{\cdot}}$, q.e.d. \hfill$\square$
 
 \medskip
 {\bf 4.6.} Take any $D_\cdot \! = (D_0 \subset \ldots \subset D_n )\in 
    \CN er_{n} T^{*}X^a$ and 
    $\{ Z_{\alpha },V_{\alpha}\}\in\CA_{D_{\cdot}}$. Let $A' $ be the subset 
 of the set $A$ of indices formed by those $\alpha$  that $Z_\alpha$ lies in $U\bar{_0}$ if $D_0$ is a lens and in $U_0$ otherwise; set  $A'':=A\smallsetminus A'$. Set 
 $K^\dagger (\{ Z_{\alpha}\}, D_0 ):=  (\mathop\prod\limits_{\alpha\in A'} 
 K^\dagger (Z_{\alpha}, D_0 ))\times 
 (\mathop\prod\limits_{\alpha\in A''} K(Z_{\alpha})) $, where 
 $K^\dagger (Z_{\alpha}, D_0)$ is $K^{\CN_0 }(Z_{\alpha})$ from 3.4 if 
 $D_{0}$ is a lens, and is the contractible $K$-spectrum $K^{triv}(Z_{\alpha})$ of
 the DG category $\CP (Z_{\alpha})^{triv}$ of acyclic complexes from $\CP (Z_{\alpha})$ otherwise. 
We leave it to the reader to check that the morphism 
$K^\dagger (\{ Z_{\alpha}\} ,D_0 )\buildrel{\Sigma i_{\alpha !}}
\over\longrightarrow 
K^{\dagger}(D_{0})$ is a weak equivalence of spectra (cf.~3.4). By 2.5, $R\Gamma_{c}$ vanishes on
each $\CP (Z_{\alpha })^{\CN_{0}}$, hence it maps $K^\dagger (Z_{\alpha},D_{0})$ to 
$K^{triv}$.

Consider two chains of morphisms of spectra
 $$ K (X)
\buildrel{\Sigma i_{\alpha !}}\over\longleftarrow 
K(\{ Z_{\alpha}\})
\buildrel{(R\Gamma_c )}\over\lra \mathop\prod\limits_{\alpha \in A} K 
 \buildrel{(tr)}\over\longleftarrow 
 \mathop\prod\limits_{\alpha \in A} C_\n (V_\alpha ,K ) \to\tag 4.6.1
$$$$\to (\mathop\prod\limits_{\alpha \in A'} 
C_\n (V_\alpha ,K ))\times (\mathop\prod\limits_{\alpha \in A''} 
C_\n (V_\alpha ,V_\alpha \smallsetminus U_{0};K ))\to K^!_X (U_{0}),$$
  $$ K^{\dagger}(D_{0})
\buildrel{\Sigma i_{\alpha !}}\over\longleftarrow 
K^\dagger (\{ Z_{\alpha}\},D_{0})
\buildrel{(R\Gamma_c )}\over\lra (\mathop\prod\limits_{\alpha \in A'} K^{triv} )
\times (\mathop\prod\limits_{\alpha \in A''} K )
 \buildrel{(tr)}\over\longleftarrow \tag 4.6.2$$ $$ \leftarrow
(\mathop\prod\limits_{\alpha \in A'} C_\n (V_\alpha ,K^{triv} ))
\times (\mathop\prod\limits_{\alpha \in A''} C_\n (V_\alpha ,K )) \to  
$$$$\to (\mathop\prod\limits_{\alpha \in A'}
C_\n (V_\alpha ,K^{triv} ))\times (\mathop\prod\limits_{\alpha \in A''} 
C_\n (V_\alpha ,V_\alpha \smallsetminus U_{0};K ))\to K^!_X (U_{0}),$$
 whose left looking arrows are weak equivalences. There is an evident morphism of diagrams
 $\iota :$ (4.6.2)$\to$(4.6.1) which is identity map at the last term. Notice that the composition of 
 (4.6.1) (with the left looking arrows homotopically inverted) is the same as that
 of $K(X)\buildrel{\varepsilon_{\{Z_{\alpha},V_{\alpha}\}}}\over\lra 
 C_{\n}(X,K)\to K^{!}_{X}(U_{0})$ (see (4.1.2)). Thus
 the composition    $\eta_{\{Z_{\alpha},V_{\alpha}\} D_\cdot}$   
  of (4.6.2) is the same as that of  $K^\dagger (D_{0})\buildrel{\iota}\over\lra K(X)\buildrel{\varepsilon_{\{Z_{\alpha},V_{\alpha}\}}}\over\lra 
 C_{\n}(X,K)\to K^{!}_{X}(U_{0})$. Denote by   $G_{\{ Z_\alpha ,V_\alpha \} D_\cdot }$ the penultimate term in (4.6.2); notice that it is homotopically trivial. 
 
 As in (4.1.3), for $\{ Z_{\beta}, V_{\beta}\}\ge \{ Z_{\alpha },V_{\alpha}\}$ in $\CA_{D_{\cdot}}$
 there are natural morphisms of diagrams $\kappa :$ (4.6.1)$_{\{ Z_{\beta},V_{\beta}\}}\to$ 
 (4.6.1)$_{\{ Z_{\alpha},V_{\alpha}\}}$,  (4.6.2)$_{\{ Z_{\beta},V_{\beta}\}}\to$
 (4.6.2)$_{\{ Z_{\alpha},V_{\alpha}\}}$ which are identity maps
 at the beginning and end terms. They commute with 
 the $\iota$ maps.
   Consider the homotopy $\CA_{D_{\cdot}}$-limit of diagrams (4.6.1) 
 and (4.6.2), and add to them the arrows
  $K(X)\to$ holim$_{\CA_{D_{\cdot}}}K(X)$,
 $K^{\dagger}(D_{0})\to$ holim$_{\CA_{D_{\cdot}}}K^\dagger (D_{0})$ from
 the left and holim$_{\CA_{D_{\cdot}}}K^{!}_{X}(U_{0})\leftarrow K^{!}_{X}(U_{0})$
  from the right (which are weak equivalences by 4.5). Let $\varepsilon_{D_{\cdot}}: K(X)\to
 K^{!}_{X}(U_{0})$,
  $\eta_{D_{\cdot}}: K^{\dagger}(D_{0})\to
 K^{!}_{X}(U_{0})$ be the compositions. One has $\eta_{D_{\cdot}}=\varepsilon_{D_{\cdot
 }}\iota$, and this morphism factors canonically through the homotopically trivial spectrum $G_{D_\cdot }:= $holim$_{\CA_{D_{\cdot}}} G_{\{ Z_\alpha ,V_\alpha \} D_\cdot }$.
 
 The above constructions are compatible with the simplicial structure on the nerve. Namely, 
if $\phi : [0,m]\to [0,n]$ is a monotone map, then $D_{\cdot}\in\CN er_{n}T^{*}X^{a}$
 yields $D^{\phi}_{\cdot}\in\CN er_{m}T^{*}X^{a}$. One has $\CA_{D_{\cdot}}\subset\CA_{D^{\phi}_{\cdot}}$,
 and for each $\{ Z_{\alpha},V_{\alpha}\}\in\CA_{D_{\cdot}}$ 
 there are evident morphisms of diagrams $\phi :$ (4.6.1)$_{
 \{ Z_{\alpha},V_{\alpha}\}D^{\phi}_{\cdot}}\to$ 
 (4.6.1)$_{\{ Z_{\alpha},V_{\alpha}\}D_{\cdot}}$,  (4.6.2)$_{
 \{ Z_{\alpha},V_{\alpha}\}D^{\phi}_{\cdot}}\to$
 (4.6.2)$_{\{ Z_{\alpha},V_{\alpha}\}D_{\cdot}}$, which commute with the $\iota$ and $\kappa$ maps.
 Passing to holim$_{\CA_{D_\cdot}}$, we get cohomological type coefficient systems on $\CN er \, T^{* }X^{a}$ (see 1.3) and  morphisms of those. 
 
Any cohomological type coefficient system $L=L_{D_\cdot}$ on $ \CN er T^{*}X^{a}$ yields a presheaf $L^a$ on $T^* X$, $L^a (E ):= C^\n (\CN er E^a ,L)$. If $L$ comes from a presheaf on $T^* X$ (as in Example (ii) in 1.3), then the corresponding $L^a$ equals that from 4.3. 

Applying this functor to $\varepsilon_{D_\cdot}$, we get the morphism $\varepsilon^a : K(X)_{T^* X}^a \to (\pi^*_{T^*} K^!_X )^a$. Applying it to  
 $\eta_{D_{\cdot}  }$, we get $\varepsilon^a \iota^a$, which is the morphism $\eta^{ a}$ from 4.4. Therefore $\eta^a$ factors through 
the  homotopically trivial $G^a$,  and we are done.
\hfill$\square$

\end